\documentclass[12pt,a4paper]{amsart}
\usepackage[all]{xy}
\usepackage{amscd}
\usepackage{amsfonts}
\usepackage{amsmath}
\usepackage{amssymb}
\usepackage{amsthm}
\usepackage{bbm}
\usepackage[colorlinks,linkcolor=blue]{hyperref}
\usepackage{calc}
\usepackage{extarrows}
\usepackage{graphicx}
\usepackage{graphics}
\usepackage{mathrsfs}
\usepackage[margin=2.9cm]{geometry}

\theoremstyle{plain}
\newtheorem{thm}{Theorem}[section]

\newtheorem{ppst}[thm]{Proposition}

\theoremstyle{definition}
\newtheorem{defn}[thm]{Definition}
\newtheorem{rem}[thm]{Remark}
\newtheorem{exam}[thm]{Example}

\title[]{Finite Hyperfields of order $n\leq 5$}
\author{Ziqi Liu}
\address{Jilin University, Changchun, Jilin, China}
\email{liuzq0616@mails.jlu.edu.cn}

\thanks{This work was started in February 2020 as an entertainment during quarantine time. The first draft was finished in April 2020. I thank Dr. Mansour Eyvazi for pointing out several mistakes in the first version of this work.}

\begin{document}

\begin{abstract}
  In this paper, the author introduces the concept and basic properties of finite (commutative) hyperfields. Also, the author shows that, up to isomorphism, there are exactly 2 hyperfields of order 2; 5 hyperfields of order 3; 7 hyperfields of order 4; 27 hyperfields of order 5. Those hyperfields could be first hand materials for those who are doing relevant researches.
\end{abstract}

\maketitle
\tableofcontents

\newpage
\section{Backgrounds}
\subsection{Hypergroups, Hypergroups and Hyperfields}
\begin{defn}
A \textbf{hyperoperation} on a set $S$ is a map $\boxplus:S\times S\rightarrow 2^S\backslash\{\varnothing\}$. Moreover, for a hyperoperation $\boxplus$ on $S$ and non-empty subsets $A,B$ of $S$, we define
$$A\boxplus B:=\bigcup_{a\in A,b\in B}(a\boxplus b)$$
A hyperoperation $\boxplus$ is called \textbf{commutative} if $a\boxplus b=b\boxplus a$ for all $a,b\in S$. If not especially mentioned, hyperoperations in this paper will always be commutative.\\
A hyperoperation $\boxplus$ is called \textbf{associative} if $a\boxplus(b\boxplus c)=(a\boxplus b)\boxplus c$ for all $a,b,c\in S$.
\end{defn}
\begin{defn}
Given an associative hyperoperation $\boxplus$, the \textbf{hypersum} of $x_1,x_2,\dots,x_m$ for $m\geq2$ is recursively defined as
$$x_1\boxplus\cdots\boxplus x_n:=\bigcup_{x'\in x_2\boxplus\cdots\boxplus x_n}x_1\boxplus x'$$
\end{defn}
\begin{defn}
A \textbf{hypergroup} is a tuple $G:=(G,\boxplus,0)$, where $\boxplus$ is an associative hyperoperation on $G$ such that:\\
(1) $0\boxplus x=\{x\}$ for all $x\in G$;\\
(2) For every $x\in G$ there is a unique element $-x$ of $G$ such that $0\in x\boxplus-x$;\\
(3) $x\in y\boxplus z$ if and only if $z\in x\boxplus(-y)$.\\
Here $-x$ is often called as the \textbf{hyperinverse} of $x$ and (3) as the reversibility axiom.
\end{defn}
\begin{defn}
A (Krasner) \textbf{hyperring} is a tuple $R:=(R,\odot,\boxplus,1,0)$ such that:\\
(1) $(R,\odot,1)$ is a commutative monoid;\\
(2) $(R,\boxplus,0)$ is a commutative hypergroup;\\
(3) $0\odot x=x\odot0=0$ for all $x\in R$;\\
(4) $a\odot(x\boxplus y)=(a\odot x)\boxplus(a\odot y)$ for all $a,x,y\in R$;\\
(5) $(x\boxplus y)\odot a=(x\odot a)\boxplus(y\odot a)$ for all $a,x,y\in R$.\\
In the following part, we often use the underlying set $R$ to refer to a hyperring and may omit $\odot$ if there is no likehood of confusion.
\end{defn}
\begin{defn}
A hyperring $F:=(F,\odot,\boxplus,1,0)$ is called a \textbf{hyperfield} if $0\neq 1$ and every non-zero element of $F$ has a multiplicative inverse.
\end{defn}

\begin{exam}
If $\mathbb{F}:=(\mathbb{F},\cdot,+,1,0)$ is a field, then $\mathbb{F}$ is isomorphic to a hyperfield $(\mathbb{F},\odot,\boxplus,1,0)$ where $x\odot y=x\cdot y$ and $x\boxplus y=\{x+y\}$ for all $x,y\in\mathbb{F}$.\\
In the following context, the term field $\mathbb{F}$ refers to the hyperfield isomorphic to $\mathbb{F}$.
\end{exam}
\begin{exam}
Consider tuple $(\{0,1\},\odot,\boxplus,1,0)$ where $\odot$ is the usual multiplication and hyperaddition $\boxplus$ is defined by
$$0\boxplus0=\{0\},\qquad1\boxplus0=0\boxplus1=\{1\},\qquad 1\boxplus1=\{0,1\}$$
then $\mathbb{K}:=(\{0,1\},\odot,\boxplus,1,0)$ is a hyperfield, called the \textbf{Krasner hyperfield}.
\end{exam}
\begin{rem}
In general, let $\Gamma=(\Gamma,\cdot,1,\leq)$ be a totally ordered abelian group, one can define a hyperfield structure $(\Gamma\cup\{0\},\odot,\boxplus,1,0)$ where $\boxplus$ is defined as
$$x\boxplus y=\left\{
\begin{aligned}
\{\max\{x,y\}\} &, & x\neq y \\
\{z\in\Gamma:z\leq x\}\cup\{0\} &,  & x=y
\end{aligned}
\right.$$
for all $x,y\in\Gamma\cup\{0\}$ and $\odot$ is an extension of $\cdot$ with $0\odot x=0$ for all $x\in\Gamma\cup\{0\}$.\\
Such a hyperfield is called a \textbf{valuative hyperfield}.
\end{rem}

\newpage
\begin{exam}
Given an abelian group $G:=(G,\cdot,1)$ and a self-inverse element $e\in G$, one can check that $W(G,e):=(G\cup\{0\},\odot,\boxplus,1,0)$ is a hyperfield where the multiplication $\odot$ is an extension of $\cdot$ in $G$, in other words, $\odot$ is generated by
$$x\odot y=\left\{
\begin{aligned}
x\cdot y &,\quad x, y\in G \\
0 &,\quad   x=0,y\in G\cup\{0\}
\end{aligned}
\right.$$
and the hyperaddition $\boxplus$ is defined by
$$0\boxplus x=\{x\},\quad x\boxplus(e\cdot x)=G\cup\{0\},\quad x\boxplus y=G$$
for any nonzero $x$ and $y$ with $y\neq ex$. $W(G,e)$ is called a \textbf{weak hyperfield}.
\end{exam}

\subsection{Homomorphisms and Isomorphisms of Hyperfields}
\begin{defn}
Given two hyperfields $\mathbb{F},\mathbb{L}$ and a map $f:\mathbb{F}\rightarrow\mathbb{L}$, if for all $x,y\in \mathbb{F}$\\
(1) $f(0)=0,f(1)=1$;\\
(2) $f(xy)=f(x)f(y)$ and $f(x\boxplus_{\mathbb{F}} y)\subseteq f(x)\boxplus_{\mathbb{L}} f(y)$;\\
then $f:\mathbb{F}\rightarrow\mathbb{L}$ is a (weak) \textbf{homomorphism of hyperfields}.
\end{defn}
\begin{exam}\label{RK}
Consider the field of real numbers $\mathbb{R}$ and the Krasner hyperfield $\mathbb{K}$, one has a natural homomorphism of hyperfields
\begin{center}
$f:\mathbb{K}\longrightarrow\mathbb{K};\quad 0\longmapsto 0,\,\,x\longmapsto 1$
\end{center}
for all $x\in\mathbb{R}-\{0\}$.
\end{exam}
\begin{defn}
Given two hyperfields $\mathbb{F},\mathbb{L}$ and a map $f:\mathbb{F}\rightarrow\mathbb{L}$, if for all $x,y\in \mathbb{F}$\\
(1) $f(0)=0,f(1)=1$;\\
(2) $f(xy)=f(x)f(y)$ and $f(x)\boxplus_{\mathbb{L}} f(y)=\{f(x\boxplus_{\mathbb{F}} y)\}$;\\
then the map $f:\mathbb{F}\rightarrow\mathbb{L}$ is a strong \textbf{homomorphism of hyperfields}.
\end{defn}
\begin{defn}
Given two hyperfields $\mathbb{F},\mathbb{L}$ and a homomorphism of hyperfields $f:\mathbb{F}\rightarrow\mathbb{L}$, if a homomorphism of hyperfields $g:\mathbb{L}\rightarrow\mathbb{F}$ satisfies:\\
(1) $f\circ g=\textup{id}_{\mathbb{L}}$;\\
(2) $g\circ f=\textup{id}_{\mathbb{F}}$;\\
then the map $f:\mathbb{F}\rightarrow\mathbb{L}$ is an \textbf{isomorphism of hyperfields} and hyperfield $\mathbb{F}$ is \textbf{isomorphic} to hyperfield $\mathbb{L}$. One can check that an isomorphism of hyperfields is a strong homomorphism of hyperfields.
\end{defn}
\begin{ppst} Given two hyperfields $\mathbb{F},\mathbb{L}$, then an isomorphism of hyperfields $f:\mathbb{F}\rightarrow\mathbb{L}$ is a bijection and a strong homomorphism of hyperfields.

\begin{proof}
By definition, there exists a homomorphism of hyperfields $g:\mathbb{L}\rightarrow\mathbb{F}$ with $f\circ g=\textup{id}_{\mathbb{L}}$, so $f$ must be surjective to cover all elements in $\mathbb{L}$. Similarly, $f$ should be injective since $g\circ f=\textup{id}_{\mathbb{F}}$ is injective. Therefore, $f$ is a bijection.\\
Since $f(x\boxplus_{\mathbb{F}} y)\subseteq f(x)\boxplus_{\mathbb{L}} f(y)$ and
$$f(x\boxplus_{\mathbb{F}}y)=f(g(f(x))\boxplus_{\mathbb{F}}g(f(y)))\supseteq f(g(f(x)\boxplus_{\mathbb{L}}f(y)))=f(x)\boxplus_{\mathbb{L}}f(y)$$
One has $f(x\boxplus_{\mathbb{F}} y)= f(x)\boxplus_{\mathbb{L}} f(y)$ for all $x,y\in \mathbb{F}$.
\end{proof}
\end{ppst}

\begin{exam}\label{F2K}
Consider the bijection
$$f:\mathbb{F}_2\rightarrow\mathbb{K};\quad x\longmapsto x$$
between the finite field $\mathbb{F}_2$ and the Krasner field $\mathbb{K}$. One is able to check that $f$ is a hyperfield homomorphism but not an isomorphism. This example shows a fact that a bijective hyperfield homomorphism is not necessary a isomorphism.
\end{exam}

\subsection{The Underlying Monoid of a Finite Hyperfield} The underlying monoid of a finite hyperfield is a crucible way to study the structure of a finite hyperfield, since what exactly such a underlying monoid should be is known.

\begin{defn}
For a hyperfield $\mathbb{F}=(F,\odot,\boxplus,1,0)$, if the underlying set $F$ is a finite set, then $\mathbb{F}$ is a \textbf{finite hyperfield}. The cardinal number of $F$ is called the \textbf{order} of $\mathbb{F}$, denoted by $|\mathbb{F}|$. The monoid $(\mathbb{F})^\times:=(F,\odot,1)$ is called the \textbf{underlying monoid} of hyperfield $\mathbb{F}$.
\end{defn}
\begin{exam} It is clear to see the following statements.\\
(1) There are no finite fields of order $1$ by our definition.\\
(2) The Karnser hyperfield $\mathbb{K}$ is a finite hyperfield of order $2$.\\
(3) All finite field $\mathbb{F}_n$ are finite hyperfields.\\
(4) For each $n>1$, there exist at least a weak hyperfield of order $n$.
\end{exam}

\begin{ppst}\label{lemma0}
For a finite hyperfield $\mathbb{F}$, if $1\neq-1$, then\\
(1) $|\mathbb{F}|$ is an odd number;\\
(2) $a\boxplus(-a)$ is a set of odd cardinal number for all $a\in\mathbb{F}$.

\begin{proof}
(1) One has $-x=x\odot(-1)\neq x\odot 1=x$ for all nonzero $x\in\mathbb{F}$ because $-1\neq1$. Therefore, $|\mathbb{F}|$ is an odd number.\\
(2) For every element $x\in a\boxplus(-a)$, one has $-x\in a\boxplus(-a)$ since
$$a\boxplus(-a)=(-a)\boxplus a=-1\odot(a\boxplus(-a))$$
Note that $0\in a\boxplus(-a)$ for all $a\in\mathbb{F}$ and $-x\neq x$ for all nonzero element $x$, it is clear that $|a\boxplus(-a)|$ is an odd number.
\end{proof}
\end{ppst}

Note that the underlying monoid of finite hyperfield can be written as $G\cup\{0\}$ where $G$ is a finite abelian group, so one needs the following theorem to help us figure out all possible underlying monoids for finite hyperfields.
\begin{thm}\label{FTA} (Fundamental Theorem of Finite Abelian Groups)
Every finite abelian group is a direct product of cyclic groups of a prime power order. Moreover,
the number of terms in the product and the orders of the cyclic groups are uniquely
determined by the group.
\end{thm}
\begin{exam}
Given $\mathbb{Z}_n:=\mathbb{Z}/n\mathbb{Z}$, then one knows that all abelian group of order 27 are isomorphic to one of the above
$$\mathbb{Z}_{27},\quad\mathbb{Z}_{9}\times\mathbb{Z}_{3},\quad\mathbb{Z}_{3}\times\mathbb{Z}_{3}\times\mathbb{Z}_{3}$$
according to Theorem \ref{FTA}.

In the following part, the group $C_n=((a),\cdot,1)$ is defined as a multiplicatively-written cyclic group of order $n$, the group $C_{i_1,\dots,i_k}$ is given by
$$C_{i_1,\dots,i_k}:=C_{i_1}\times\cdots\times C_{i_k}$$
It is clear that $C_n\cong\mathbb{Z}_n$ and $C_{i_1,\dots,i_k}\cong\mathbb{Z}_{i_1}\times\cdots\times\mathbb{Z}_{i_k}$. Moreover, for hyperfields with the same underlying monoid, we use their hyperaddtions to represent them.
\end{exam}

\subsection{Hyperfield Extensions} Analogy of the concept field extension, we introduce the concept of hyperfield extension and discuss some special features for hyperfield extension of finite hyperfields.

\begin{defn}
Given two hyperfields $\mathbb{F}$ and $\mathbb{L}$, if there exists an injective strong hyperfield homomorphism $f:\mathbb{F}\rightarrow\mathbb{L}$, then hyperfield $\mathbb{L}$ is called a (weak) \textbf{hyperfield extension} of $\mathbb{F}$, denoted $\mathbb{L}/\mathbb{F}$.
\end{defn}
\begin{exam}
Consider the weak hyperfield $W(C_2,1)$ where $C_2:=(\{1,a\},\cdot)$ is the cyclic group of order 2, the homomorphism of hyperfields $f:\mathbb{K}\rightarrow W(C_2,1)$ given by $f(x)=x$ is clearly injective. Then $W(C_2,1)$ is an extension of $\mathbb{K}$.
\end{exam}
\begin{exam}
The homomorphism of hyperfields described in Example \ref{F2K} implies that the Krasner field $\mathbb{K}$ is an extension of finite field $\mathbb{F}_2$.
\end{exam}
\begin{rem} As one can define different hyperfield structures based on monoid $C_n\cup\{0\}$, an identity map can be a hyperfield homomorphism of different hyperfields.
\end{rem}

\begin{defn}
Given two hyperfields $\mathbb{F}$ and $\mathbb{L}$, if there exists an injective strong hyperfield homomorphism $f:\mathbb{F}\rightarrow\mathbb{L}$, then hyperfield extension $\mathbb{L}/\mathbb{F}$ is a \textbf{strong hyperfield extension}.
\end{defn}
\begin{exam}
Finite field $\mathbb{F}_{2^2}$ is a strong hyperfield extension of finite field $\mathbb{F}_2$. In fact, all homomorphisms of finite fields are strong hyperfield homomorphisms.
\end{exam}

\section{Finite Hyperfields of order $n\leq 5$}In this part, we will give all hyperfields of order $n$ and show the number of different hyperfields of order $n$ up to isomorphism for $n\leq 5$.
\subsection{Finite Hyperfields of order 2}
\begin{ppst}\label{F2}
Finite field $\mathbb{F}_2$ and the Krasner hyperfield $\mathbb{K}$ are the only two hyperfields of order $2$.

\begin{proof}
It is clear that the underlying set of a hyperfield $\mathbb{F}$ with order $2$ must be $\{0,1\}$ and it follows that $1=-1$. Then consider the set $1\boxplus1$, since $0$ is an element of it by definition, one has either $1\boxplus 1=\{0\}$ or $1\boxplus 1=\{0,1\}$. The first one gives $\mathbb{F}_2$ and the second one gives $\mathbb{K}$.
\end{proof}
\end{ppst}
\subsection{Finite Hyperfields of order 3}
\begin{defn}
Given an abelian group $G=(\{1,-1\},\cdot,1)$, define a binary operation $\odot$ as the extension of $\cdot$ on $G\cup\{0\}$ and a hyperaddition $\boxplus$ on $G\cup\{0\}$ by
$$\begin{tabular}{|c|c|c|c|}
  \hline
  $\boxplus$ & $0$ & $1$ & $-1$\\
  \hline
  $0$ & $\{0\}$ & $\{1\}$ & $\{-1\}$\\
  \hline
  $1$ & $\{1\}$ & $\{1\}$ & $\{0,1,-1\}$\\
  \hline
  $-1$ & $\{-1\}$ & $\{0,1,-1\}$ & $\{-1\}$\\
\hline
\end{tabular}$$
Then $\mathbb{S}:=(\{0,1,-1\},\odot,\boxplus,1,0)$ is a hyperfield, called the \textbf{hyperfield of signs}.
\end{defn}
\begin{defn}
Given the cyclic group $C_2:=(\{1,a\},\cdot,1)$, define a special weak hyperfield $\mathbb{W}:=W(C_2,-1)$ to be \textbf{weak hyperfield of signs}.
\end{defn}

\begin{ppst}\label{F3}
There are exactly five hyperfields of order $3$.

\begin{proof}
The underlying monoid of a hyperfield of order $3$ must be $C_2\cup\{0\}$.\\
\textbf{Case 1} When $a$ is the hyperinverse of $1$, we write the underlying set as $\{1,-1,0\}$. By Proposition \ref{lemma0}, set $1\boxplus(-1)$ is $\{0\}$ or $\{1,-1,0\}$.\\
\textbf{Case 1.1}  If $1\boxplus(-1)=\{1,-1,0\}$, we only have $\mathbb{S}$ and $\mathbb{W}$.\\
\textbf{Case 1.2} If $1\boxplus(-1)=\{0\}$, then $1\notin1\boxplus1$ and then $1\boxplus1=\{-1\}$. So, the only possible hyperfield here is finite field $\mathbb{F}_3$. Note that $\mathbb{S}$ is not an extension of $\mathbb{F}_3$ and
$$\mathbb{F}_3\longrightarrow \mathbb{W} \longleftarrow \mathbb{S}$$
\textbf{Case 2} When $1$ is self-inverse, we write the underlying set as $\{1,a,0\}$. Since in this case $1\in1\boxplus a\Leftrightarrow a\in1\boxplus a$, $1\boxplus a=\{1,a\}$. So we only need to consider the set $1\boxplus1$.\\
\textbf{Case 2.1} If $1\boxplus1=\{0\}$, we have $1\notin1\boxplus a$ since $a\notin 1\boxplus1$, a contradiction.\\
\textbf{Case 2.2} If $1\boxplus1=\{a,0\}$, one can check that the hyperaddition defined as below
$$\begin{tabular}{|c|c|c|c|}
  \hline
  $\boxplus$ & $0$ & $1$ & $a$\\
  \hline
  $0$ & $\{0\}$ & $\{1\}$ & $\{a\}$\\
  \hline
  $1$ & $\{1\}$ & $\{a,0\}$ & $\{1,a\}$\\
  \hline
  $a$ & $\{a\}$ & $\{1,a\}$ & $\{1,0\}$\\
\hline
\end{tabular}$$
gives a hyperfield structure on $C_2\cup\{0\}$. This hyperfield can be seen as a natural extension of $\mathbb{F}_2$ by adding an element $a$.
$$\begin{tabular}{|c|c|c|c|}
  \hline
  $\boxplus$ & $0$ & $1$ & $\quad$\\
  \hline
  $0$ & $\{0\}$ & $\{1\}$ &  $\quad$\\
  \hline
  $1$ & $\{1\}$ & $\{0\}$ &  $\quad$\\
  \hline
   & & & $\quad$ \\
\hline
\end{tabular}\longrightarrow
\begin{tabular}{|c|c|c|c|}
  \hline
  $\boxplus$ & $0$ & $1$ & $a$\\
  \hline
  $0$ & $\{0\}$ & $\{1\}$ & $\{a\}$\\
  \hline
  $1$ & $\{1\}$ & $\{a,0\}$ & $\{1,a\}$\\
  \hline
  $a$ & $\{a\}$ & $\{1,a\}$ & $\{1,0\}$\\
\hline
\end{tabular}$$
We denote this hyperfield as $(\mathbb{F}_2)^{\uparrow3}$.\\
\textbf{Case 2.3} If $1\boxplus1=\{1,a,0\}$, one can check that the hyperaddition defined as below
$$\begin{tabular}{|c|c|c|c|}
  \hline
  $\boxplus$ & $0$ & $1$ & $a$\\
  \hline
  $0$ & $\{0\}$ & $\{1\}$ & $\{a\}$\\
  \hline
  $1$ & $\{1\}$ & $\{1,a,0\}$ & $\{1,a\}$\\
  \hline
  $a$ & $\{a\}$ & $\{1,a\}$ & $\{1,a,0\}$\\
\hline
\end{tabular}$$
gives a hyperfield structure on $C_2\cup\{0\}$, which is actually $W(C_2,1)$. This hyperfield can be seen as a natural extension of $\mathbb{F}_2$ by adding an element $a$.
$$\begin{tabular}{|c|c|c|c|}
  \hline
  $\boxplus$ & $0$ & $1$ & $\quad$\\
  \hline
  $0$ & $\{0\}$ & $\{1\}$ &  $\quad$\\
  \hline
  $1$ & $\{1\}$ & $\{1,0\}$ &  $\quad$\\
  \hline
   & & & $\quad$ \\
\hline
\end{tabular}\longrightarrow
\begin{tabular}{|c|c|c|c|}
  \hline
  $\boxplus$ & $0$ & $1$ & $a$\\
  \hline
  $0$ & $\{0\}$ & $\{1\}$ & $\{a\}$\\
  \hline
  $1$ & $\{1\}$ & $\{1,a,0\}$ & $\{1,a\}$\\
  \hline
  $a$ & $\{a\}$ & $\{1,a\}$ & $\{1,a,0\}$\\
\hline
\end{tabular}$$
We denote this hyperfield as $\mathbb{K}^{\uparrow3}:=W(C_2,1)$.
\end{proof}
\end{ppst}
\begin{rem}
We have the following commutative diagram
$$
\xymatrix{
\mathbb{F}_2 \ar[d]_{\textup{natural extension by adding }a} \ar[rr]^{\mathfrak{i}_2}&& \mathbb{K} \ar[d]^{\textup{natural extension by adding }a} \\
(\mathbb{F}_2)^{\uparrow3} \ar[rr]^{\mathfrak{i}_3} &&  \mathbb{K}^{\uparrow3}
}
$$
where $\mathfrak{i}_2$ and $\mathfrak{i}_3$ are identity map on set $C_1\cup\{0\}$ and $C_2\cup\{0\}$ respectively.
\end{rem}

\subsection{Finite Hyperfields of order 4}
The underlying monoid of a hyperfield of order $4$ should be $C_3\cup\{0\}=\{0,1,a,a^2\}$. According to Proposition \ref{lemma0}, $1=-1$ in all finite hyperfield of order $4$. Here we should pay more attention to the associativity property of a hyperaddition in a hyperfield since it from now on will make big differences on the numbers of different finite hyperfield of a certain order.

\begin{ppst}\label{F4}
There are exactly seven hyperfields of order $4$.

\begin{proof}
We divide cases by the set $1\boxplus1$.\\
\textbf{Case 1} If $1\boxplus 1=\{0\}$, then $1\boxplus a=\{a^2\}$ since $\{1,a\}\cap 1\boxplus a=\varnothing$, which implies that the only possible hyperaddition in this case is
$$\begin{tabular}{|c|c|c|c|c|}
  \hline
  $\boxplus$ & $0$ & $1$ & $a$ & $a^2$\\
  \hline
  $0$ & $\{0\}$ & $\{1\}$ & $\{a\}$ & $\{a^2\}$\\
  \hline
  $1$ & $\{1\}$ & $\{0\}$ & $\{a^2\}$ & $\{a\}$\\
  \hline
  $a$ & $\{a\}$ & $\{a^2\}$ & $\{0\}$ & $\{1\}$\\
  \hline
  $a^2$ & $\{a^2\}$ & $\{a\}$ & $\{1\}$ & $\{0\}$\\
\hline
\end{tabular}$$
which gives finite field $\mathbb{F}_{2^2}$.\\
\textbf{Case 2} If $1\boxplus 1=\{1,0\}$, then $1\boxplus a=\{a^2\}$ since $\{1,a\}\cap 1\boxplus a=\varnothing$, which implies that $1\boxplus a^2=\{a\}$. However, in this case $1\boxplus(a\boxplus1)=1\boxplus a^2=\{a\}$ but
$$a\boxplus(1\boxplus1)=a\boxplus \{1,0\}=\{a,a^2\}$$
a contradiction. Therefore, there are on valid hyperfields in this case.\\
\textbf{Case 3} If $1\boxplus 1=\{a,0\}$, then $1\in1\boxplus a$ and $a\notin 1\boxplus a$.\\
\textbf{Case 3.1} When $1\boxplus a=\{1\}$, $1\boxplus(1\boxplus a)=1\boxplus1=\{a,0\}$ and
$$a\boxplus(1\boxplus1)=a\boxplus\{a,0\}=\{a,a^2,0\}$$
which is impossible.\\
\textbf{Case 3.2} When $1\boxplus a=\{1,a^2\}$, we have
$$\begin{tabular}{|c|c|c|c|c|}
  \hline
  $\boxplus$ & $0$ & $1$ & $a$ & $a^2$\\
  \hline
  $0$ & $\{0\}$ & $\{1\}$ & $\{a\}$ & $\{a^2\}$\\
  \hline
  $1$ & $\{1\}$ & $\{a,0\}$ & $\{1,a^2\}$ & $\{a,a^2\}$\\
  \hline
  $a$ & $\{a\}$ & $\{1,a^2\}$ & $\{a^2,0\}$ & $\{1,a\}$\\
  \hline
  $a^2$ & $\{a^2\}$ & $\{a,a^2\}$ & $\{1,a\}$ & $\{1,0\}$\\
\hline
\end{tabular}$$
\textbf{Case 4} If $1\boxplus 1=\{a^2,0\}$, then $a\in1\boxplus a$ and $1\notin 1\boxplus a$.\\
\textbf{Case 4.1} When $1\boxplus a=\{a\}$, $1\boxplus(1\boxplus a)=1\boxplus a=\{a\}$ and
$$a\boxplus(1\boxplus1)=a\boxplus\{a^2,0\}=a(1\boxplus a)\cup\{a\}=\{a,a^2\}$$
which is impossible.\\
\textbf{Case 4.2} When $1\boxplus a=\{a,a^2\}$, we have
$$\begin{tabular}{|c|c|c|c|c|}
  \hline
  $\boxplus$ & $0$ & $1$ & $a$ & $a^2$\\
  \hline
  $0$ & $\{0\}$ & $\{1\}$ & $\{a\}$ & $\{a^2\}$\\
  \hline
  $1$ & $\{1\}$ & $\{a^2,0\}$ & $\{a,a^2\}$ & $\{1,a\}$\\
  \hline
  $a$ & $\{a\}$ & $\{a,a^2\}$ & $\{1,0\}$ & $\{1,a^2\}$\\
  \hline
  $a^2$ & $\{a^2\}$ & $\{1,a\}$ & $\{1,a^2\}$ & $\{a,0\}$\\
\hline
\end{tabular}$$
Here we should notice that the two hyperfields given in \textbf{Case 3.2} and \textbf{Case 4.2} are isomorphically the same one, since it is easy to check that the map
$$f:C_3\cup\{0\}\rightarrow C_3\cup\{0\},\quad (1,a,a^2,0)\mapsto(1,a^2,a,0)$$
is the isomorphism between them. Its inverse is the map itself. Moreover, this hyperfield is an extension of finite hyperfield $\mathbb{F}_2$, we denote it by $(\mathbb{F}_2)^{\uparrow 4}$.\\
\textbf{Case 5} If $1\boxplus 1=\{1,a,0\}$, then $1\in1\boxplus a$ and $a\notin 1\boxplus a$. Here we have
$$a\boxplus(1\boxplus1)=a\boxplus\{1,a,0\}=1\boxplus a\cup\{a,a^2,0\}=\{1,a,a^2,0\}$$
therefore $a^2\in 1\boxplus a$ or $1\boxplus(1\boxplus a)=1\boxplus1\neq\{1,a,a^2,0\}$, a contradiction. In this case, $1\boxplus a=\{1,a^2\}$ and then we have
$$\begin{tabular}{|c|c|c|c|c|}
  \hline
  $\boxplus$ & $0$ & $1$ & $a$ & $a^2$\\
  \hline
  $0$ & $\{0\}$ & $\{1\}$ & $\{a\}$ & $\{a^2\}$\\
  \hline
  $1$ & $\{1\}$ & $\{1,a,0\}$ & $\{1,a^2\}$ & $\{a,a^2\}$\\
  \hline
  $a$ & $\{a\}$ & $\{1,a^2\}$ & $\{a,a^2,0\}$ & $\{1,a\}$\\
  \hline
  $a^2$ & $\{a^2\}$ & $\{a,a^2\}$ & $\{1,a\}$ & $\{1,a^2,0\}$\\
\hline
\end{tabular}$$
\textbf{Case 6} If $1\boxplus 1=\{1,a^2,0\}$, then $a\in1\boxplus a$ and $1\notin 1\boxplus a$. Here we have
$$a\boxplus(1\boxplus1)=a\boxplus\{1,a^2,0\}\supseteq\{a,a^2\}$$
therefore $a^2\in 1\boxplus a$ or $1\boxplus(1\boxplus a)=1\boxplus a=\{a\}\neq a\boxplus(1\boxplus1)$, a contradiction. In this case, $1\boxplus a=\{a,a^2\}$ and then we have
$$\begin{tabular}{|c|c|c|c|c|}
  \hline
  $\boxplus$ & $0$ & $1$ & $a$ & $a^2$\\
  \hline
  $0$ & $\{0\}$ & $\{1\}$ & $\{a\}$ & $\{a^2\}$\\
  \hline
  $1$ & $\{1\}$ & $\{1,a^2,0\}$ & $\{a,a^2\}$ & $\{1,a\}$\\
  \hline
  $a$ & $\{a\}$ & $\{a,a^2\}$ & $\{1,a,0\}$ & $\{1,a^2\}$\\
  \hline
  $a^2$ & $\{a^2\}$ & $\{1,a\}$ & $\{1,a^2\}$ & $\{a,a^2,0\}$\\
\hline
\end{tabular}$$
Here we should notice that the two hyperfields given in \textbf{Case 5} and \textbf{Case 6} are isomorphically the same one, since it is easy to check that the map
$$f:C_3\cup\{0\}\rightarrow C_3\cup\{0\},\quad (1,a,a^2,0)\mapsto(1,a^2,a,0)$$
is the isomorphism between them. Its inverse is the map itself. Moreover, this hyperfield is an extension of the Krasner hyperfield $\mathbb{K}$, we denote it by $\mathbb{K}^{\uparrow 4}$.\\
\textbf{Case 7} If $1\boxplus 1=\{a,a^2,0\}$, then $\{1,a\}\subseteq1\boxplus a$. In this case, we could obtain the following two hyperadditions:
$$\begin{tabular}{|c|c|c|c|c|}
  \hline
  $\boxplus$ & $0$ & $1$ & $a$ & $a^2$\\
  \hline
  $0$ & $\{0\}$ & $\{1\}$ & $\{a\}$ & $\{a^2\}$\\
  \hline
  $1$ & $\{1\}$ & $\{a,a^2,0\}$ & $\{1,a\}$ & $\{1,a^2\}$\\
  \hline
  $a$ & $\{a\}$ & $\{1,a\}$ & $\{1,a^2,0\}$ & $\{a,a^2\}$\\
  \hline
  $a^2$ & $\{a^2\}$ & $\{1,a^2\}$ & $\{a,a^2\}$ & $\{1,a,0\}$\\
\hline
\end{tabular}$$
and
$$\begin{tabular}{|c|c|c|c|c|}
  \hline
  $\boxplus$ & $0$ & $1$ & $a$ & $a^2$\\
  \hline
  $0$ & $\{0\}$ & $\{1\}$ & $\{a\}$ & $\{a^2\}$\\
  \hline
  $1$ & $\{1\}$ & $\{a,a^2,0\}$ & $\{1,a,a^2\}$ & $\{1,a,a^2\}$\\
  \hline
  $a$ & $\{a\}$ & $\{1,a,a^2\}$ & $\{1,a^2,0\}$ & $\{1,a,a^2\}$\\
  \hline
  $a^2$ & $\{a^2\}$ & $\{1,a,a^2\}$ & $\{1,a,a^2\}$ & $\{1,a,0\}$\\
\hline
\end{tabular}$$
These two hyperfields are extensions of $(\mathbb{F}_2)^{\uparrow4}$, and we denote them by $(\mathbb{F}_2)^{\uparrow^24}$ and $(\mathbb{F}_2)^{\uparrow^24}_{\rightarrow}$ respectively. Note that we also have $(\mathbb{F}_2)^{\uparrow^24}_{\rightarrow}/(\mathbb{F}_2)^{\uparrow^24}$.\\
\textbf{Case 8} If $1\boxplus 1=\{1,a,a^2,0\}$, then $\{1,a\}\subseteq1\boxplus a$. In this case, we could obtain the following two hyperadditions:
$$\begin{tabular}{|c|c|c|c|c|}
  \hline
  $\boxplus$ & $0$ & $1$ & $a$ & $a^2$\\
  \hline
  $0$ & $\{0\}$ & $\{1\}$ & $\{a\}$ & $\{a^2\}$\\
  \hline
  $1$ & $\{1\}$ & $\{1,a,a^2,0\}$ & $\{1,a\}$ & $\{1,a^2\}$\\
  \hline
  $a$ & $\{a\}$ & $\{1,a\}$ & $\{1,a,a^2,0\}$ & $\{a,a^2\}$\\
  \hline
  $a^2$ & $\{a^2\}$ & $\{1,a^2\}$ & $\{a,a^2\}$ & $\{1,a,a^2,0\}$\\
\hline
\end{tabular}$$
and a weak hyperfield
$$\begin{tabular}{|c|c|c|c|c|}
  \hline
  $\boxplus$ & $0$ & $1$ & $a$ & $a^2$\\
  \hline
  $0$ & $\{0\}$ & $\{1\}$ & $\{a\}$ & $\{a^2\}$\\
  \hline
  $1$ & $\{1\}$ & $\{1,a,a^2,0\}$ & $\{1,a,a^2\}$ & $\{1,a,a^2\}$\\
  \hline
  $a$ & $\{a\}$ & $\{1,a,a^2\}$ & $\{1,a,a^2,0\}$ & $\{1,a,a^2\}$\\
  \hline
  $a^2$ & $\{a^2\}$ & $\{1,a,a^2\}$ & $\{1,a,a^2\}$ & $\{1,a,a^2,0\}$\\
\hline
\end{tabular}$$
These two hyperfields are extensions of $\mathbb{K}$, so we denote them by $\mathbb{K}^{\uparrow^24}$ and $\mathbb{K}^{\uparrow^24}_{\rightarrow}$ respectively. Note that $\mathbb{K}^{\uparrow^24}_{\rightarrow}/\mathbb{K}^{\uparrow^24}$ and $\mathbb{K}^{\uparrow^24}_{\rightarrow}$ is also the weak hyperfield $W(C_3,1)$.
\end{proof}
\end{ppst}
\begin{rem}
We have the following commutative diagram
$$
\xymatrix{
  &&  \mathbb{F}_{2^2} \ar[r] & (\mathbb{F}_2)^{\uparrow4} \ar[ddd] \ar[dr] &\\
\mathbb{F}_2 \ar[d] \ar[urr] \ar[rr] && (\mathbb{F}_2)^{\uparrow^24} \ar[d] \ar[rr] && (\mathbb{F}_2)^{\uparrow^24}_{\rightarrow} \ar[d] \\
\mathbb{K} \ar[rr]  && \mathbb{K}^{\uparrow^24} \ar[rr] && \mathbb{K}^{\uparrow^24}_{\rightarrow}\\
   &&&  \mathbb{K}^{\uparrow4} \ar[ur] &
}
$$
One can also observe that all hyperfields of order 4 are extensions of known hyperfields.
\end{rem}

\subsection{Finite Hyperfields of order 5}
The underlying monoid of a hyperfield of order $5$ should be $C_4\cup\{0\}$ or $C_{2,2}\cup\{0\}$. When the underlying monoid is $C_4\cup\{0\}$, 1 and $a^2$ are only two possible options for the inverse element of $1$. When the underlying monoid is $C_{2,2}\cup\{0\}:=\{1,a,b,ab,0\}$, only $0$ cannot be the inverse of $1$. Without losing generality, we only need to consider $1,ab$ as the inverse element of $1$.

\begin{ppst}\label{F5,1}
There are nine different hyperfields with $a^2=-1$ of underlying monoid $C_4\cup\{0\}$.

\begin{proof}
According to Proposition \ref{lemma0}, $1\boxplus a^2$ should be one of the following sets
$$\{0\},\,\,\{0,1,a^2\},\,\,\{0,a,a^3\},\,\,\{0,1,a,a^2,a^3\}$$
\textbf{Case 1} If $1\boxplus a^2=\{0\}$, one can see
\begin{center}
$\{a,a^3\}\cap 1\boxplus a^2=\varnothing$ $\Leftrightarrow$ $1\notin1\boxplus a$ and $a\notin 1\boxplus a$
\end{center}
Therefore, $1\boxplus a$ might be $\{a^2\},\{a^3\}$ or $\{a^2,a^3\}$.\\
\textbf{Case 1.1} When $1\boxplus a$ is $\{a^2\}$ or $\{a^3\}$, we could obtain exactly one hyperfield isomorphic to finite field $\mathbb{F}_5$.\\
\textbf{Case 1.2} When $1\boxplus a=\{a^2,a^3\}$, $\{a,a^3\}\subseteq1\boxplus1$ and $a^3\in a\boxplus a^2$. Therefore, we have
$$a^3\in a^2\boxplus \{a,a^3\}\subseteq a^2\boxplus(1\boxplus 1) $$
But $1\boxplus(1\boxplus a^2)=1\boxplus 0=\{1\}$, so we have no valid hyperfields in this case.\\
\textbf{Case 2} If $1\boxplus a^2=\{0,1,a^2\}$, then $1\boxplus a$ also should be one of $\{a^2\},\{a^3\}$ and $\{a^2,a^3\}$.\\
\textbf{Case 2.1} When $1\boxplus a=\{a^2\}$, $\{1,a^3\}\subseteq 1\boxplus1\subseteq\{1,a^2,a^3\}$. In this case, we have
$$a\boxplus(1\boxplus 1)\supseteq a\boxplus\{1,a^3\}=\{0,a,a^2,a^3\}$$
but $1\boxplus(1\boxplus a)=1\boxplus\{a^2\}=\{0,1,a^2\}$, a contradiction.

\newpage
\noindent\textbf{Case 2.2} When $1\boxplus a=\{a^3\}$, $\{1,a\}\subseteq 1\boxplus1\subseteq\{1,a,a^2\}$. In this case, we have
$$a^3\in a\boxplus\{1,a\}\subseteq  a\boxplus(1\boxplus 1)=1\boxplus(1\boxplus a)= 1\boxplus\{a^3\}$$
but $a^3\in 1\boxplus a^3$ is equivalent to $a^3\in a^2\boxplus a^3$ and then $a\in 1\boxplus a$, a contradiction.\\
\textbf{Case 2.3} When $1\boxplus a=\{a^2,a^3\}$, $\{1,a,a^3\}\subseteq 1\boxplus1\subseteq\{1,a,a^2,a^3\}$. We have
$$a^3\in a\boxplus\{1,a,a^3\}\subseteq  a\boxplus(1\boxplus 1)=1\boxplus(1\boxplus a)= 1\boxplus\{a^2,a^3\}=\{0,1,a,a^2\}$$
which is impossible. Hence there are no valid hyperfields in this case and then, as a conclusion, in the whole \textbf{Case 2}.\\
\textbf{Case 3} If $1\boxplus a^2=\{0,a,a^3\}$, then $\{1,a\}\subseteq 1\boxplus a$. We divide cases by $1\boxplus a$.\\
\textbf{Case 3.1} When $1\boxplus a=\{1,a\}$, we have $\{a,a^3\}\cap1\boxplus 1=\varnothing$. Notice that $1\notin1\boxplus1$ since $1\notin1\boxplus a^2$, it must be $1\boxplus 1=\{a^2\}$. One can see $(1\boxplus1)\boxplus a^2=a^2\boxplus a^2=1$, but
$$a\in(1\boxplus a^2)\boxplus 1=\{0,a,a^3\}\boxplus1=a^3\boxplus 1\cup\{1,a\}$$
a contradiction.\\
\textbf{Case 3.2} When $1\boxplus a=\{1,a,a^2\}$, we have $\{a^3\}\subseteq1\boxplus 1\subseteq\{a^2,a^3\}$. Moreover, one can see $a\in1\boxplus a^3\Leftrightarrow a^2\in1\boxplus a$ and $a^2\notin1\boxplus a^3\Leftrightarrow a^3\notin1\boxplus a$. However, we have
$$(1\boxplus1)\boxplus a\subseteq\{a^2,a^3\}\boxplus a=\{1,a,a^2,0\}$$
while
$$a^2\in(1\boxplus a)\boxplus 1=\{1,a,a^2\}\boxplus 1=1\boxplus1\cap\{a,a^2\}\boxplus1$$
\textbf{Case 3.3} When $1\boxplus a=\{1,a,a^3\}$, one can get the following two hyperfields
$$\begin{tabular}{|c|c|c|c|c|c|}
  \hline
  $\boxplus$ & $0$ & $1$ & $a$ & $a^2$ & $a^3$\\
  \hline
  $0$ & $\{0\}$ & $\{1\}$ & $\{a\}$ & $\{a^2\}$ & $\{a^3\}$\\
  \hline
  $1$ & $\{1\}$ & $\{a,a^2\}$ & $\{1,a,a^3\}$ & $\{0,a,a^3\}$ & $\{1,a^2,a^3\}$\\
  \hline
  $a$ & $\{a\}$ & $\{1,a,a^3\}$ & $\{a^2,a^3\}$ & $\{1,a,a^2\}$ & $\{0,1,a^2\}$\\
  \hline
  $a^2$ & $\{a^2\}$ & $\{0,a,a^3\}$ & $\{1,a,a^2\}$ & $\{1,a^3\}$ & $\{a,a^2,a^3\}$\\
  \hline
  $a^3$ & $\{a^3\}$ & $\{1,a^2,a^3\}$ & $\{0,1,a^2\}$ & $\{a,a^2,a^3\}$ & $\{1,a\}$\\
  \hline
\end{tabular}$$
and
$$\begin{tabular}{|c|c|c|c|c|c|}
  \hline
  $\boxplus$ & $0$ & $1$ & $a$ & $a^2$ & $a^3$\\
  \hline
  $0$ & $\{0\}$ & $\{1\}$ & $\{a\}$ & $\{a^2\}$ & $\{a^3\}$\\
  \hline
  $1$ & $\{1\}$ & $\{a^2,a^3\}$ & $\{1,a,a^2\}$ & $\{0,a,a^3\}$ & $\{1,a,a^3\}$\\
  \hline
  $a$ & $\{a\}$ & $\{1,a,a^2\}$ & $\{1,a^3\}$ & $\{a,a^2,a^3\}$ & $\{0,1,a^2\}$\\
  \hline
  $a^2$ & $\{a^2\}$ & $\{0,a,a^3\}$ & $\{a,a^2,a^3\}$ & $\{1,a\}$ & $\{1,a^2,a^3\}$\\
  \hline
  $a^3$ & $\{a^3\}$ & $\{1,a,a^3\}$ & $\{0,1,a^2\}$ & $\{1,a^2,a^3\}$ & $\{a,a^2\}$\\
  \hline
\end{tabular}$$
These two hyperfield are isomorphic, we denote it by $\mathbb{Y}$.\\
\textbf{Case 3.4} When $1\boxplus a=\{1,a,a^2,a^3\}$, we have
$$\{a,a^3\}\subseteq1\boxplus 1\subseteq\{a,a^2,a^3\}$$
Moreover, one can see $1\boxplus a^3=\{1,a,a^2,a^3\}$. So if $1\boxplus1=\{a,a^3\}$, we have
$$a\boxplus(1\boxplus1)=1\boxplus\{a,a^3\}=\{1,a^2,0\}$$
while
$$a\in (1\boxplus a)\boxplus 1=\{1,a,a^2,a^3\}\boxplus 1=\{1,a,a^2,a^3,0\}$$
a contradiction.

\noindent Therefore, $1\boxplus 1=\{a,a^2,a^3\}$ and we have
$$\begin{tabular}{|c|c|c|c|c|c|}
  \hline
  $\boxplus$ & $0$ & $1$ & $a$ & $a^2$ & $a^3$\\
  \hline
  $0$ & $\{0\}$ & $\{1\}$ & $\{a\}$ & $\{a^2\}$ & $\{a^3\}$\\
  \hline
  $1$ & $\{1\}$ & $\{a,a^2,a^3\}$ & $\{1,a,a^2,a^3\}$ & $\{0,a,a^3\}$ & $\{1,a,a^2,a^3\}$\\
  \hline
  $a$ & $\{a\}$ & $\{1,a,a^2,a^3\}$ & $\{1,a^2,a^3\}$ & $\{1,a,a^2,a^3\}$ & $\{0,1,a^2\}$\\
  \hline
  $a^2$ & $\{a^2\}$ & $\{0,a,a^3\}$ & $\{1,a,a^2,a^3\}$ & $\{1,a,a^3\}$ & $\{1,a,a^2,a^3\}$\\
  \hline
  $a^3$ & $\{a^3\}$ & $\{1,a,a^2,a^3\}$ & $\{0,1,a^2\}$ & $\{1,a,a^2,a^3\}$ & $\{1,a,a^2\}$\\
  \hline
\end{tabular}$$
This hyperfield is an extension of $\mathbb{Y}$, denoted by $\mathbb{Y}_\rightarrow$.\\
\textbf{Case 4} If $1\boxplus a^2=\{0,1,a,a^2,a^3\}$, then $\{1,a\}\subseteq 1\boxplus a$. We divide cases by $1\boxplus a$.\\
\textbf{Case 4.1} When $1\boxplus a=\{1,a\}$, $\{a,a^3\}\cap1\boxplus 1=\varnothing$. Note that $1\in 1\boxplus1\Leftrightarrow 1\in 1\boxplus a^2$, we could obtain the following two hyperfields.
$$\begin{tabular}{|c|c|c|c|c|c|}
  \hline
  $\boxplus$ & $0$ & $1$ & $a$ & $a^2$ & $a^3$\\
  \hline
  $0$ & $\{0\}$ & $\{1\}$ & $\{a\}$ & $\{a^2\}$ & $\{a^3\}$\\
  \hline
  $1$ & $\{1\}$ & $\{1\}$ & $\{1,a\}$ & $\{0,1,a,a^2,a^3\}$ & $\{1,a^3\}$\\
  \hline
  $a$ & $\{a\}$ & $\{1,a\}$ & $\{a\}$ & $\{a,a^2\}$ & $\{0,1,a,a^2,a^3\}$\\
  \hline
  $a^2$ & $\{a^2\}$ & $\{0,1,a,a^2,a^3\}$ & $\{a,a^2\}$ & $\{a^2\}$ & $\{a^2,a^3\}$\\
  \hline
  $a^3$ & $\{a^3\}$ & $\{1,a^3\}$ & $\{0,1,a,a^2,a^3\}$ & $\{a^2,a^3\}$ & $\{a^3\}$\\
  \hline
\end{tabular}$$
and
$$\begin{tabular}{|c|c|c|c|c|c|}
  \hline
  $\boxplus$ & $0$ & $1$ & $a$ & $a^2$ & $a^3$\\
  \hline
  $0$ & $\{0\}$ & $\{1\}$ & $\{a\}$ & $\{a^2\}$ & $\{a^3\}$\\
  \hline
  $1$ & $\{1\}$ & $\{1,a^2\}$ & $\{1,a\}$ & $\{0,1,a,a^2,a^3\}$ & $\{1,a^3\}$\\
  \hline
  $a$ & $\{a\}$ & $\{1,a\}$ & $\{a,a^3\}$ & $\{a,a^2\}$ & $\{0,1,a,a^2,a^3\}$\\
  \hline
  $a^2$ & $\{a^2\}$ & $\{0,1,a,a^2,a^3\}$ & $\{a,a^2\}$ & $\{1,a^2\}$ & $\{a^2,a^3\}$\\
  \hline
  $a^3$ & $\{a^3\}$ & $\{1,a^3\}$ & $\{0,1,a,a^2,a^3\}$ & $\{a^2,a^3\}$ & $\{a,a^3\}$\\
  \hline
\end{tabular}$$
The first hyperfield is an extension of the hyperfield of signs $\mathbb{S}$, we denote it by $\mathbb{S}^{\uparrow5}$. The second hyperfield is an extension of the weak hyperfield of signs $\mathbb{W}$, we denote it by $\mathbb{W}^{\uparrow5}$. The hyperfield homomorphism is $f:(0,1,-1)\mapsto(0,1,a^2)$.\\
\textbf{Case 4.2} When $1\boxplus a=\{1,a,a^2\}$, we have $\{1,a^3\}\subseteq1\boxplus 1\subseteq\{1,a^2,a^3\}$. Moreover, one can see $a\in1\boxplus a^3\Leftrightarrow a^2\in1\boxplus a$ and $a^2\notin1\boxplus a^3\Leftrightarrow a^3\notin1\boxplus a$. In this case, we could obtain the following hyperfields.
$$\begin{tabular}{|c|c|c|c|c|c|}
  \hline
  $\boxplus$ & $0$ & $1$ & $a$ & $a^2$ & $a^3$\\
  \hline
  $0$ & $\{0\}$ & $\{1\}$ & $\{a\}$ & $\{a^2\}$ & $\{a^3\}$\\
  \hline
  $1$ & $\{1\}$ & $\{1,a^3\}$ & $\{1,a,a^2\}$ & $\{0,1,a,a^2,a^3\}$ & $\{1,a,a^3\}$\\
  \hline
  $a$ & $\{a\}$ & $\{1,a,a^2\}$ & $\{1,a\}$ & $\{a,a^2,a^3\}$ & $\{0,1,a,a^2,a^3\}$\\
  \hline
  $a^2$ & $\{a^2\}$ & $\{0,1,a,a^2,a^3\}$ & $\{a,a^2,a^3\}$ & $\{a,a^2\}$ & $\{1,a^2,a^3\}$\\
  \hline
  $a^3$ & $\{a^3\}$ & $\{1,a,a^3\}$ & $\{0,1,a,a^2,a^3\}$ & $\{1,a^2,a^3\}$ & $\{a^2,a^3\}$\\
  \hline
\end{tabular}$$
and
$$\begin{tabular}{|c|c|c|c|c|c|}
  \hline
  $\boxplus$ & $0$ & $1$ & $a$ & $a^2$ & $a^3$\\
  \hline
  $0$ & $\{0\}$ & $\{1\}$ & $\{a\}$ & $\{a^2\}$ & $\{a^3\}$\\
  \hline
  $1$ & $\{1\}$ & $\{1,a^2,a^3\}$ & $\{1,a,a^2\}$ & $\{0,1,a,a^2,a^3\}$ & $\{1,a,a^3\}$\\
  \hline
  $a$ & $\{a\}$ & $\{1,a,a^2\}$ & $\{1,a,a^3\}$ & $\{a,a^2,a^3\}$ & $\{0,1,a,a^2,a^3\}$\\
  \hline
  $a^2$ & $\{a^2\}$ & $\{0,1,a,a^2,a^3\}$ & $\{a,a^2,a^3\}$ & $\{1,a,a^2\}$ & $\{1,a^2,a^3\}$\\
  \hline
  $a^3$ & $\{a^3\}$ & $\{1,a,a^3\}$ & $\{0,1,a,a^2,a^3\}$ & $\{1,a^2,a^3\}$ & $\{a,a^2,a^3\}$\\
  \hline
\end{tabular}$$
We denote these two hyperfields by $\mathbb{S}^{\uparrow5}_\rightarrow$ and $\mathbb{W}^{\uparrow5}_\rightarrow$ respectively.

\newpage
\noindent\textbf{Case 4.3} When $1\boxplus a=\{1,a,a^3\}$, we have $\{1,a\}\subseteq1\boxplus 1\subseteq\{1,a,a^2\}$. Moreover, one can see $a^2\in1\boxplus a^3\Leftrightarrow a^3\in1\boxplus a$ and $a\notin1\boxplus a^3\Leftrightarrow a^2\notin1\boxplus a$. In this case, we could obtain the following hyperfields.
$$\begin{tabular}{|c|c|c|c|c|c|}
  \hline
  $\boxplus$ & $0$ & $1$ & $a$ & $a^2$ & $a^3$\\
  \hline
  $0$ & $\{0\}$ & $\{1\}$ & $\{a\}$ & $\{a^2\}$ & $\{a^3\}$\\
  \hline
  $1$ & $\{1\}$ & $\{1,a\}$ & $\{1,a,a^3\}$ & $\{0,1,a,a^2,a^3\}$ & $\{1,a^2,a^3\}$\\
  \hline
  $a$ & $\{a\}$ & $\{1,a,a^3\}$ & $\{a,a^2\}$ & $\{1,a,a^2\}$ & $\{0,1,a,a^2,a^3\}$\\
  \hline
  $a^2$ & $\{a^2\}$ & $\{0,1,a,a^2,a^3\}$ & $\{1,a,a^2\}$ & $\{a^2,a^3\}$ & $\{a,a^2,a^3\}$\\
  \hline
  $a^3$ & $\{a^3\}$ & $\{1,a^2,a^3\}$ & $\{0,1,a,a^2,a^3\}$ & $\{a,a^2,a^3\}$ & $\{1,a^3\}$\\
  \hline
\end{tabular}$$
and
$$\begin{tabular}{|c|c|c|c|c|c|}
  \hline
  $\boxplus$ & $0$ & $1$ & $a$ & $a^2$ & $a^3$\\
  \hline
  $0$ & $\{0\}$ & $\{1\}$ & $\{a\}$ & $\{a^2\}$ & $\{a^3\}$\\
  \hline
  $1$ & $\{1\}$ & $\{1,a,a^2\}$ & $\{1,a,a^3\}$ & $\{0,1,a,a^2,a^3\}$ & $\{1,a^2,a^3\}$\\
  \hline
  $a$ & $\{a\}$ & $\{1,a,a^3\}$ & $\{a,a^2,a^3\}$ & $\{1,a,a^2\}$ & $\{0,1,a,a^2,a^3\}$\\
  \hline
  $a^2$ & $\{a^2\}$ & $\{0,1,a,a^2,a^3\}$ & $\{1,a,a^2\}$ & $\{1,a^2,a^3\}$ & $\{a,a^2,a^3\}$\\
  \hline
  $a^3$ & $\{a^3\}$ & $\{1,a^2,a^3\}$ & $\{0,1,a,a^2,a^3\}$ & $\{a,a^2,a^3\}$ & $\{1,a,a^3\}$\\
  \hline
\end{tabular}$$
One can see clearly that these hyperfields are isomorphic to $\mathbb{S}^{\uparrow5}_\rightarrow$ and $\mathbb{W}^{\uparrow5}_\rightarrow$ respectively.\\
\textbf{Case 4.4} When $1\boxplus a=\{1,a,a^2,a^3\}$, we have
$$\{1,a,a^3\}\subseteq1\boxplus 1\subseteq\{1,a,a^2,a^3\}$$
Moreover, one can see $1\boxplus a^3=\{1,a,a^2,a^3\}$. Then we could obtain two hyperfields.
$$\begin{tabular}{|c|c|c|c|c|c|}
  \hline
  $\boxplus$ & $0$ & $1$ & $a$ & $a^2$ & $a^3$\\
  \hline
  $0$ & $\{0\}$ & $\{1\}$ & $\{a\}$ & $\{a^2\}$ & $\{a^3\}$\\
  \hline
  $1$ & $\{1\}$ & $\{1,a,a^3\}$ & $\{1,a,a^2,a^3\}$ & $\{0,1,a,a^2,a^3\}$ & $\{1,a,a^2,a^3\}$\\
  \hline
  $a$ & $\{a\}$ & $\{1,a,a^2,a^3\}$ & $\{1,a,a^2\}$ & $\{1,a,a^2,a^3\}$ & $\{0,1,a,a^2,a^3\}$\\
  \hline
  $a^2$ & $\{a^2\}$ & $\{0,1,a,a^2,a^3\}$ & $\{1,a,a^2,a^3\}$ & $\{a,a^2,a^3\}$ & $\{1,a,a^2,a^3\}$\\
  \hline
  $a^3$ & $\{a^3\}$ & $\{1,a,a^2,a^3\}$ & $\{0,1,a,a^2,a^3\}$ & $\{1,a,a^2,a^3\}$ & $\{1,a^2,a^3\}$\\
  \hline
\end{tabular}$$
and
$$\begin{tabular}{|c|c|c|c|c|c|}
  \hline
  $\boxplus$ & $0$ & $1$ & $a$ & $a^2$ & $a^3$\\
  \hline
  $0$ & $\{0\}$ & $\{1\}$ & $\{a\}$ & $\{a^2\}$ & $\{a^3\}$\\
  \hline
  $1$ & $\{1\}$ & $\{1,a,a^2,a^3\}$ & $\{1,a,a^2,a^3\}$ & $\{0,1,a,a^2,a^3\}$ & $\{1,a,a^2,a^3\}$\\
  \hline
  $a$ & $\{a\}$ & $\{1,a,a^2,a^3\}$ & $\{1,a,a^2,a^3\}$ & $\{1,a,a^2,a^3\}$ & $\{0,1,a,a^2,a^3\}$\\
  \hline
  $a^2$ & $\{a^2\}$ & $\{0,1,a,a^2,a^3\}$ & $\{1,a,a^2,a^3\}$ & $\{1,a,a^2,a^3\}$ & $\{1,a,a^2,a^3\}$\\
  \hline
  $a^3$ & $\{a^3\}$ & $\{1,a,a^2,a^3\}$ & $\{0,1,a,a^2,a^3\}$ & $\{1,a,a^2,a^3\}$ & $\{1,a,a^2,a^3\}$\\
  \hline
\end{tabular}$$
The first one is an extension of $\mathbb{S}^{\uparrow5}_\rightarrow$ and the second one is an extension of $\mathbb{W}^{\uparrow5}_\rightarrow$. We denote them by $\mathbb{S}^{\uparrow5}_{\rightarrow^2}$ and $\mathbb{W}^{\uparrow5}_{\rightarrow^2}$ respectively.
\end{proof}
\end{ppst}
\begin{rem}
We have the following commutative diagram
$$
\xymatrix{
 &&\mathbb{S} \ar[rr] \ar[d] && \mathbb{S}^{\uparrow5} \ar[rr] \ar[d] && \mathbb{S}^{\uparrow5}_\rightarrow \ar[rr] \ar[d] && \mathbb{S}^{\uparrow5}_{\rightarrow^2} \ar[d]\\
\mathbb{F}_3\ar[rr]&&\mathbb{W} \ar[rr] && \mathbb{W}^{\uparrow5} \ar[rr] && \mathbb{W}^{\uparrow5}_\rightarrow \ar[rr] && \mathbb{W}^{\uparrow5}_{\rightarrow^2}\\
&&&&\mathbb{F}_5\ar[uurr]&&\mathbb{Y}\ar[rr]&&\mathbb{Y}_\rightarrow\ar[u]
}
$$
\end{rem}

\begin{ppst}\label{F5,2}
There are seven different hyperfields with $1=-1$ of underlying monoid $C_4\cup\{0\}$.

\begin{proof}
Here we divide cases by the set $1\boxplus1$ and focus on associativity.\\
\textbf{Case 1} If $1\boxplus 1=\{0\}$, then $\{1,a\}\cap1\boxplus a=\varnothing$. Note that
$$a^2\in1\boxplus a\Leftrightarrow 1\in a\boxplus a^2\Leftrightarrow a^3\in 1\boxplus a$$
So $1\boxplus a=\{a^2,a^3\}$, then $1\boxplus a^3=\{a,a^2\}$ and $1\boxplus a^2=\{a,a^3\}$. Now we have
$$1\boxplus(1\boxplus a)=1\boxplus\{a^2,a^3\}=\{a,a^2,a^3\}$$
but $a\boxplus(1\boxplus1)=a\boxplus0=\{a\}$, a contradiction.\\
\textbf{Case 2} If $1\boxplus 1=\{1,0\}$, then $1\boxplus a=\{a^2,a^3\}$ still holds. Therefore we have
\begin{center}$1\boxplus a^3=\{a,a^2\}$ and $1\boxplus a^2=\{a,a^3\}$\end{center}
In this case, we could obtain a hyperfield as following.
$$\begin{tabular}{|c|c|c|c|c|c|}
 \hline
  $\boxplus$ & $0$ & $1$ & $a$ & $a^2$ & $a^3$\\
  \hline
  $0$ & $\{0\}$ & $\{1\}$ & $\{a\}$ & $\{a^2\}$ & $\{a^3\}$\\
  \hline
  $1$ & $\{1\}$ & $\{1,0\}$ & $\{a^2,a^3\}$ & $\{a,a^3\}$ & $\{a,a^2\}$\\
  \hline
  $a$ & $\{a\}$  & $\{a^2,a^3\}$ & $\{a,0\}$ & $\{1,a^3\}$ & $\{1,a^2\}$\\
  \hline
  $a^2$ & $\{a^2\}$ & $\{a,a^3\}$ & $\{1,a^3\}$ & $\{a^2,0\}$ & $\{1,a\}$\\
  \hline
  $a^3$ & $\{a^3\}$ & $\{a,a^2\}$ & $\{1,a^2\}$ & $\{1,a\}$ & $\{a^3,0\}$\\
  \hline
\end{tabular}$$
This hyperfield is an extension of $\mathbb{K}$, denoted $\mathbb{K}^{\uparrow5}$.\\
\textbf{Case 3} If $1\boxplus 1=\{a,0\}$, then $\{1,a^2\}\cap1\boxplus a^2=\varnothing$. Therefore, $1\boxplus a^2=\{a,a^3\}$ and $1\boxplus a=\{1,a^2,a^3\}$ or otherwise $1\boxplus a^2$ would be empty. Hence
$$1\boxplus (1\boxplus a)=1\boxplus\{1,a^2,a^3\}=\{a,a^2,a^3,0\}$$
while $a\boxplus(1\boxplus 1)=a\boxplus\{a,0\}=\{a,a^2,0\}$, a contradiction.\\
\textbf{Case 4} If $1\boxplus 1=\{a^2,0\}$, then $\{1,a\}\cap1\boxplus a=\varnothing$, which implies that $1\boxplus a=\{a^2,a^3\}$. Therefore, we know $a\boxplus (1\boxplus 1)=a\boxplus\{a^2,0\}=\{1,a,a^3\}$. However, one can see
$$a^2\in 1\boxplus (1\boxplus a)=1\boxplus\{a^2,a^3\}$$
because $a^2\in 1\boxplus a^2$. Hence there are no hyperfields in this case.\\
\textbf{Case 5} If $1\boxplus 1=\{a^3,0\}$, then $1\boxplus a^2=\{a,a^3\}$ or it will be empty. In this case, we could obtain that $1\boxplus a=\{a,a^2,a^3\}$ and $1\boxplus a^3=\{1,a,a^3\}$. Then we have
$$1\boxplus (1\boxplus a^3)=1\boxplus\{1,a,a^3\}=\{a,a^2,a^3,0\}$$
and $a^3\boxplus(1\boxplus 1)=a\boxplus\{a^3,0\}=\{1,a,a^2,0\}$, a contradiction.\\
\textbf{Case 6} If $1\boxplus 1=\{1,a,0\}$, then $\{1,a^2\}\cap1\boxplus a^2=\varnothing$ and $1\boxplus a^2=\{a,a^3\}$. In this case, we have $(a\boxplus a^3)\boxplus 1=\{1,a^2\}\boxplus 1=\{1,a,a^2,0\}$ while
$$a^3\in(1\boxplus a)\boxplus a^3=\{1,a,0\}\boxplus a^3$$
a contradiction.\\
\noindent\textbf{Case 7} If $1\boxplus 1=\{1,a^2,0\}$, then $\{1,a\}\cap1\boxplus a=\varnothing$ and $1\boxplus a=\{a^2,a^3\}$. Then, we are able to obtain that $1\boxplus a^2=\{1,a,a^2,a^3\}$ and $1\boxplus a^3=\{a,a^2\}$. In this case, we have
$$(1\boxplus a^3)\boxplus a^2=\{a,a^2\}\boxplus a^2=\{1,a^2,a^3,0\}$$
while $(1\boxplus a^2)\boxplus a^3=\{1,a,a^2,a^3\}\boxplus a^3=\{1,a,a^2,a^3,0\}$
a contradiction.\\
\textbf{Case 8} If $1\boxplus 1=\{1,a^3,0\}$, then $\{1,a\}\cap1\boxplus a$ and $1\boxplus a=\{a,a^2,a^3\}$. Also, one can see $1\boxplus a^2=\{a,a^3\}$ and the find that $a^3\in(1\boxplus a)\boxplus a^3\neq(1\boxplus a^3)\boxplus a$.\\
\textbf{Case 9} If $1\boxplus 1=\{a,a^2,0\}$, then $\{1,a^2\}\subseteq1\boxplus a^2$. Here, one can see $1\boxplus a^2=C_4$ or otherwise $1\boxplus a^2=\{1,a^2\}, 1\boxplus a=\{1\}$ and hence
$$a^3\in a\boxplus(1\boxplus1)=a\boxplus\{a,a^2,0\},1\boxplus(1\boxplus a)=\{a,a^2,0\}$$
which is impossible. Therefore, $1\boxplus a^2=C_4$. We could get the following hyperfield.
$$\begin{tabular}{|c|c|c|c|c|c|}
 \hline
  $\boxplus$ & $0$ & $1$ & $a$ & $a^2$ & $a^3$\\
  \hline
  $0$ & $\{0\}$ & $\{1\}$ & $\{a\}$ & $\{a^2\}$ & $\{a^3\}$\\
  \hline
  $1$ & $\{1\}$ & $\{a,a^2,0\}$ & $\{1,a^2,a^3\}$ & $\{1,a,a^2,a^3\}$ & $\{a,a^2,a^3\}$\\
  \hline
  $a$ & $\{a\}$  & $\{1,a^2,a^3\}$ & $\{a^2,a^3,0\}$ & $\{1,a,a^3\}$ & $\{1,a,a^2,a^3\}$\\
  \hline
  $a^2$ & $\{a^2\}$ & $\{1,a,a^2,a^3\}$ & $\{1,a,a^3\}$ & $\{1,a^3,0\}$ & $\{1,a,a^2\}$\\
  \hline
  $a^3$ & $\{a^3\}$ & $\{a,a^2,a^3\}$ & $\{1,a,a^2,a^3\}$ & $\{1,a,a^2\}$ & $\{1,a,0\}$\\
  \hline
\end{tabular}$$
This is an extension of $(\mathbb{F}_2)^{\uparrow3}$, denoted $(\mathbb{F}_2)^{\uparrow^25}$.\\
\textbf{Case 10} If $1\boxplus 1=\{a,a^3,0\}$, then $1\boxplus a^2=\{a,a^3\}$ or it would be empty. In this case, one can see that $1\boxplus a=\{1,a,a^2,a^3\}$. Therefore, we have
$$(1\boxplus a^2)\boxplus a=\{a,a^3\}\boxplus a=\{1,a^2,0\}$$
while $(1\boxplus a)\boxplus a^2=\{1,a,a^2,a^3\}\boxplus a^2=\{1,a,a^2,a^3\}$, a contradiction.\\
\textbf{Case 11} If $1\boxplus 1=\{a^2,a^3,0\}$, then $1\boxplus a^2=\{1,a,a^2,a^3\}$. In this case, we have
$$\begin{tabular}{|c|c|c|c|c|c|}
 \hline
  $\boxplus$ & $0$ & $1$ & $a$ & $a^2$ & $a^3$\\
  \hline
  $0$ & $\{0\}$ & $\{1\}$ & $\{a\}$ & $\{a^2\}$ & $\{a^3\}$\\
  \hline
  $1$ & $\{1\}$ & $\{a,a^2,0\}$ & $\{1,a^2,a^3\}$ & $\{1,a,a^2,a^3\}$ & $\{a,a^2,a^3\}$\\
  \hline
  $a$ & $\{a\}$  & $\{1,a^2,a^3\}$ & $\{a^2,a^3,0\}$ & $\{1,a,a^3\}$ & $\{1,a,a^2,a^3\}$\\
  \hline
  $a^2$ & $\{a^2\}$ & $\{1,a,a^2,a^3\}$ & $\{1,a,a^3\}$ & $\{1,a^3,0\}$ & $\{1,a,a^2\}$\\
  \hline
  $a^3$ & $\{a^3\}$ & $\{a,a^2,a^3\}$ & $\{1,a,a^2,a^3\}$ & $\{1,a,a^2\}$ & $\{1,a,0\}$\\
  \hline
\end{tabular}$$
This hyperfield is isomorphic to $(\mathbb{F}_2)^{\uparrow^25}$ we mentioned before.\\
\textbf{Case 12} If $1\boxplus 1=\{1,a,a^3,0\}$, then we have $1\boxplus a^2=\{a,a^3\}$ and then have
$$(1\boxplus a^2)\boxplus a^3=\{a,a^3\}\boxplus a^3=\{1,a^2,a^3,0\}$$
while
$$a\in(1\boxplus a^3)\boxplus a^2=\{1,a,a^2,a^3\}\boxplus a^2$$
So there are no hyperfields in this case.\\
\textbf{Case 13} If $1\boxplus 1=\{1,a,a^2,0\}$, then $1\boxplus a=\{1,a^2,a^3\}$ or otherwise $1\boxplus a=\{1\}$ and one can see $$1\boxplus(a\boxplus a)\neq a\boxplus(1\boxplus a)=\{1\}$$
which will lead to contradiction. In this case, we have one hyperfield as following.
$$\begin{tabular}{|c|c|c|c|c|c|}
 \hline
  $\boxplus$ & $0$ & $1$ & $a$ & $a^2$ & $a^3$\\
  \hline
  $0$ & $\{0\}$ & $\{1\}$ & $\{a\}$ & $\{a^2\}$ & $\{a^3\}$\\
  \hline
  $1$ & $\{1\}$ & $\{1,a,a^2,0\}$ & $\{1,a^2,a^3\}$ & $\{1,a,a^2,a^3\}$ & $\{a,a^2,a^3\}$\\
  \hline
  $a$ & $\{a\}$  & $\{1,a^2,a^3\}$ & $\{a,a^2,a^3,0\}$ & $\{1,a,a^3\}$ & $\{1,a,a^2,a^3\}$\\
  \hline
  $a^2$ & $\{a^2\}$ & $\{1,a,a^2,a^3\}$ & $\{1,a,a^3\}$ & $\{1,a^2,a^3,0\}$ & $\{1,a,a^2\}$\\
  \hline
  $a^3$ & $\{a^3\}$ & $\{a,a^2,a^3\}$ & $\{1,a,a^2,a^3\}$ & $\{1,a,a^2\}$ & $\{1,a,a^3,0\}$\\
  \hline
\end{tabular}$$
This hyperfield is an extension of $(\mathbb{F}_2)^{\uparrow^25}$ and $\mathbb{K}^{\uparrow3}$, denoted $\mathbb{K}^{\uparrow^25}$.\\
\textbf{Case 14} If $1\boxplus 1=\{1,a^2,a^3,0\}$, then $1\boxplus a=\{a,a^2,a^3\}$ or otherwise $1\boxplus a=\{a\}$ and hence $a\boxplus(1\boxplus 1)\neq1\boxplus(1\boxplus a)=\{a\}$, a contradiction. In this case, we have one hyperfield as following.
$$\begin{tabular}{|c|c|c|c|c|c|}
 \hline
  $\boxplus$ & $0$ & $1$ & $a$ & $a^2$ & $a^3$\\
  \hline
  $0$ & $\{0\}$ & $\{1\}$ & $\{a\}$ & $\{a^2\}$ & $\{a^3\}$\\
  \hline
  $1$ & $\{1\}$ & $\{1,a^2,a^3,0\}$ & $\{a,a^2,a^3\}$ & $\{1,a,a^2,a^3\}$ & $\{1,a,a^2\}$\\
  \hline
  $a$ & $\{a\}$  & $\{a,a^2,a^3\}$ & $\{1,a,a^3,0\}$ & $\{1,a^2,a^3\}$ & $\{1,a,a^2,a^3\}$\\
  \hline
  $a^2$ & $\{a^2\}$ & $\{1,a,a^2,a^3\}$ & $\{1,a^2,a^3\}$ & $\{1,a,a^2,0\}$ & $\{1,a,a^3\}$\\
  \hline
  $a^3$ & $\{a^3\}$ & $\{1,a,a^2\}$ & $\{1,a,a,a^3\}$ & $\{1,a,a^3\}$ & $\{a,a^2,a^3,0\}$\\
  \hline
\end{tabular}$$
This hyperfield is isomorphic to $\mathbb{K}^{\uparrow^25}$.\\
\textbf{Case 15} If $1\boxplus 1=\{a,a^2,a^3,0\}$, then we could obtain two hyperfields. The first one is also mentioned in \cite{6} and relevant researches about quotient hyperfields.
$$\begin{tabular}{|c|c|c|c|c|c|}
 \hline
  $\boxplus$ & $0$ & $1$ & $a$ & $a^2$ & $a^3$\\
  \hline
  $0$ & $\{0\}$ & $\{1\}$ & $\{a\}$ & $\{a^2\}$ & $\{a^3\}$\\
  \hline
  $1$ & $\{1\}$ & $\{a,a^2,a^3,0\}$ & $\{1,a\}$ & $\{1,a^2\}$ & $\{1,a^3\}$\\
  \hline
  $a$ & $\{a\}$  & $\{1,a\}$ & $\{1,a^2,a^3,0\}$ & $\{a,a^2\}$ & $\{a,a^3\}$\\
  \hline
  $a^2$ & $\{a^2\}$ & $\{1,a^2\}$ & $\{a,a^2\}$ & $\{1,a,a^3,0\}$ & $\{a^2,a^3\}$\\
  \hline
  $a^3$ & $\{a^3\}$ & $\{1,a^3\}$ & $\{a,a^3\}$ & $\{a^2,a^3\}$ & $\{1,a,a^2,0\}$\\
  \hline
\end{tabular}$$
It is not isomorphic to any quotient of fields and we usually denote it by $\mathbb{M}$. Then one can see the second hyperfield in this case.
$$\begin{tabular}{|c|c|c|c|c|c|}
 \hline
  $\boxplus$ & $0$ & $1$ & $a$ & $a^2$ & $a^3$\\
  \hline
  $0$ & $\{0\}$ & $\{1\}$ & $\{a\}$ & $\{a^2\}$ & $\{a^3\}$\\
  \hline
  $1$ & $\{1\}$ & $\{a,a^2,a^3,0\}$ & $\{1,a,a^2,a^3\}$ & $\{1,a,a^2,a^3\}$ & $\{1,a,a^2,a^3\}$\\
  \hline
  $a$ & $\{a\}$  & $\{1,a,a^2,a^3\}$ & $\{1,a^2,a^3,0\}$ & $\{1,a,a^2,a^3\}$ & $\{1,a,a^2,a^3\}$\\
  \hline
  $a^2$ & $\{a^2\}$ & $\{1,a,a^2,a^3\}$ & $\{1,a,a^2,a^3\}$ & $\{1,a,a^3,0\}$ & $\{1,a,a^2,a^3\}$\\
  \hline
  $a^3$ & $\{a^3\}$ & $\{1,a,a^2,a^3\}$ & $\{1,a,a^2,a^3\}$ & $\{1,a,a^2,a^3\}$ & $\{1,a,a^2,0\}$\\
  \hline
\end{tabular}$$
It is of course an extension of $\mathbb{M}$ and also $(\mathbb{F}_2)^{\uparrow^25}$, we denote it by $(\mathbb{F}_2)^{\uparrow^25}_\rightarrow$.\\
\textbf{Case 16} If $1\boxplus 1=\{1,a,a^2,a^3,0\}$, then we could obtain the following two hyperfield.
$$\begin{tabular}{|c|c|c|c|c|c|}
 \hline
  $\boxplus$ & $0$ & $1$ & $a$ & $a^2$ & $a^3$\\
  \hline
  $0$ & $\{0\}$ & $\{1\}$ & $\{a\}$ & $\{a^2\}$ & $\{a^3\}$\\
  \hline
  $1$ & $\{1\}$ & $\{1,a,a^2,a^3,0\}$ & $\{1,a\}$ & $\{1,a^2\}$ & $\{1,a^3\}$\\
  \hline
  $a$ & $\{a\}$  & $\{1,a\}$ & $\{1,a,a^2,a^3,0\}$ & $\{a,a^2\}$ & $\{a,a^3\}$\\
  \hline
  $a^2$ & $\{a^2\}$ & $\{1,a^2\}$ & $\{a,a^2\}$ & $\{1,a,a^2,a^3,0\}$ & $\{a^2,a^3\}$\\
  \hline
  $a^3$ & $\{a^3\}$ & $\{1,a^3\}$ & $\{a,a^3\}$ & $\{a^2,a^3\}$ & $\{1,a,a^2,a^3,0\}$\\
  \hline
\end{tabular}$$
and weak hyperfield
$$\begin{tabular}{|c|c|c|c|c|c|}
 \hline
  $\boxplus$ & $0$ & $1$ & $a$ & $a^2$ & $a^3$\\
  \hline
  $0$ & $\{0\}$ & $\{1\}$ & $\{a\}$ & $\{a^2\}$ & $\{a^3\}$\\
  \hline
  $1$ & $\{1\}$ & $\{1,a,a^2,a^3,0\}$ & $\{1,a,a^2,a^3\}$ & $\{1,a,a^2,a^3\}$ & $\{1,a,a^2,a^3\}$\\
  \hline
  $a$ & $\{a\}$  & $\{1,a,a^2,a^3\}$ & $\{1,a,a^2,a^3,0\}$ & $\{1,a,a^2,a^3\}$ & $\{1,a,a^2,a^3\}$\\
  \hline
  $a^2$ & $\{a^2\}$ & $\{1,a,a^2,a^3\}$ & $\{1,a,a^2,a^3\}$ & $\{1,a,a^2,a^3,0\}$ & $\{1,a,a^2,a^3\}$\\
  \hline
  $a^3$ & $\{a^3\}$ & $\{1,a,a^2,a^3\}$ & $\{1,a,a^2,a^3\}$ & $\{1,a,a^2,a^3\}$ & $\{1,a,a^2,a^3,0\}$\\
  \hline
\end{tabular}$$
We denote the first one by $\mathbb{M}_+$ and the second one by $(\mathbb{K})^{\uparrow^25}_{\rightarrow}$.
\end{proof}
\end{ppst}
\begin{rem}
One can check the following commutative diagram.
$$
\xymatrix{
  &  & \mathbb{M}\ar[r]\ar[dr] & \mathbb{M}_+\ar[ddr] &\\
\mathbb{F}_2\ar[d]\ar[r]\ar[urr] & (\mathbb{F}_2)^{\uparrow3}\ar[d]\ar[r] & (\mathbb{F}_2)^{\uparrow^25}\ar[dr]\ar[r] & (\mathbb{F}_2)^{\uparrow^25}_\rightarrow\ar[dr] & \\
\mathbb{K}\ar[r]\ar[drr] & \mathbb{K}^{\uparrow3}\ar[rrr] && \mathbb{K}^{\uparrow^25}\ar[r]& \mathbb{K}^{\uparrow^25}_{\rightarrow}\\
 &  & \mathbb{K}^{\uparrow^5} \ar[urr] &
}
$$
Here most of the hyperfield homomorphisms are the identity map and the map send the tuple $(1,a,a^2,a^3,0)$ to $(1,a^3,a^2,a,0)$.
\end{rem}

\begin{ppst}\label{F5,3}
There are six different hyperfields with $ab=-1$ of underlying monoid $C_{2,2}\cup\{0\}$.

\begin{proof}
Here $b=(ab)a=-a$, and we divide cases by set $1\boxplus1$.\\
\textbf{Case 1} If $1\boxplus1=\{1\}$, then $\{b,ab\}\cap1\boxplus a=\varnothing$ and thus $1\boxplus a=\{1,a\}$. In this case, we have the following hyperfield.
$$\begin{tabular}{|c|c|c|c|c|c|}
 \hline
  $\boxplus$ & $0$ & $1$ & $a$ & $b$ & $ab$\\
  \hline
  $0$ & $\{0\}$ & $\{1\}$ & $\{a\}$ & $\{b\}$ & $\{ab\}$\\
  \hline
  $1$ & $\{1\}$ & $\{1\}$ & $\{1,a\}$ & $\{1,b\}$ & $\{1,a,b,ab,0\}$\\
  \hline
  $a$ & $\{a\}$  & $\{1,a\}$ & $\{a\}$ & $\{1,a,b,ab,0\}$ & $\{a,ab\}$\\
  \hline
  $b$ & $\{b\}$ & $\{1,b\}$ & $\{1,a,b,ab,0\}$ & $\{b\}$ & $\{b,ab\}$\\
  \hline
  $ab$ & $\{ab\}$ & $\{1,a,b,ab,0\}$ & $\{a,ab\}$ & $\{b,ab\}$ & $\{ab\}$\\
  \hline
\end{tabular}$$
This hyperfield is an extension of the hyperfield of sighs $\mathbb{S}$, as one can see
$$\begin{tabular}{|c|c|c|c|c|}
 \hline
  $\boxplus$ &  $1$ &  &  & $ab$\\
  \hline
  $1$ &  $\{1\}$ &  &  & $\{1,ab,0\}$\\
  \hline
   &  &  &  & \\
  \hline
   &   &  &  & \\
  \hline
  $ab$ &  $\{1,ab,0\}$ &  &  & $\{ab\}$\\
  \hline
\end{tabular}\longrightarrow \textup{this hyperfield}
$$
We denote this hyperfield by $\mathbb{S}^{\Uparrow5}$.\\
\textbf{Case 2} If $1\boxplus1=\{a\}$, then $\{b,ab\}\cap1\boxplus a=\varnothing$ and thus $1\boxplus a=\{1,a\}$. In this case,
$$1\boxplus(1\boxplus a)=1\boxplus\{1,a\}=\{1,a\}$$
but $a\boxplus(1\boxplus1)=a\boxplus a=\{1\}$, a contradiction.\\
\textbf{Case 3} If $1\boxplus1=\{b\}$, then $\{a,ab\}\cap1\boxplus b=\varnothing$ and thus $1\boxplus b=\{1,b\}$. In this case,
$$1\boxplus(1\boxplus b)=1\boxplus\{1,b\}=\{1,b\}$$
but $b\boxplus(1\boxplus1)=b\boxplus b=\{1\}$, a contradiction.\\
\textbf{Case 4} If $1\boxplus1=\{ab\}$, then $1\boxplus a=\{1,a\}$ and $1\boxplus b=\{1,b\}$. Therefore,
$$1\boxplus(1\boxplus a)=1\boxplus\{1,a\}=\{1,a,ab\}$$
but $a\boxplus(1\boxplus1)=a\boxplus ab=a(1\boxplus b)=\{a,ab\}$, a contradiction.\\
\textbf{Case 5} If $1\boxplus1=\{1,ab\}$, then $\{b,ab\}\cap1\boxplus a=\varnothing$ and thus $1\boxplus a=\{1,a\}$. In this case, we have the following hyperfield.
$$\begin{tabular}{|c|c|c|c|c|c|}
 \hline
  $\boxplus$ & $0$ & $1$ & $a$ & $b$ & $ab$\\
  \hline
  $0$ & $\{0\}$ & $\{1\}$ & $\{a\}$ & $\{b\}$ & $\{ab\}$\\
  \hline
  $1$ & $\{1\}$ & $\{1,ab\}$ & $\{1,a\}$ & $\{1,b\}$ & $\{1,a,b,ab,0\}$\\
  \hline
  $a$ & $\{a\}$  & $\{1,a\}$ & $\{a,b\}$ & $\{1,a,b,ab,0\}$ & $\{a,ab\}$\\
  \hline
  $b$ & $\{b\}$ & $\{1,b\}$ & $\{1,a,b,ab,0\}$ & $\{b,a\}$ & $\{b,ab\}$\\
  \hline
  $ab$ & $\{ab\}$ & $\{1,a,b,ab,0\}$ & $\{a,ab\}$ & $\{b,ab\}$ & $\{1,ab\}$\\
  \hline
\end{tabular}$$
This hyperfield is an extension of weak hyperfield of signs $\mathbb{W}$, denoted by $\mathbb{W}^{\Uparrow5}$.\\
\textbf{Case 6} If $1\boxplus1=\{a,ab\}$, then $\{b,ab\}\cap1\boxplus a=\varnothing$ and thus $1\boxplus a=\{1,a\}$. Hence
$$b\in a\boxplus(1\boxplus1)=a\boxplus \{a,ab\}=\{1,b\}\cup a\boxplus ab$$
but $1\boxplus(1\boxplus a)=1\boxplus\{1,a\}=\{1,a,ab\}$, a contradiction.\\
\textbf{Case 7} If $1\boxplus1=\{b,ab\}$, then $\{a,ab\}\cap1\boxplus b=\varnothing$ and thus $1\boxplus a=\{1,b\}$. Hence
$$a\in b\boxplus(1\boxplus1)=b\boxplus \{b,ab\}=\{1,a\}\cup b\boxplus ab$$
but $1\boxplus(1\boxplus b)=1\boxplus\{1,b\}=\{1,b,ab\}$, a contradiction.

\noindent\textbf{Case 8} If $1\boxplus1=\{a,b\}$, then $\{1,ab\}\cap1\boxplus ab=\varnothing$ and $1\boxplus ab=\{a,b\}$. In this case, we have $1\boxplus a=\{1,a,b,ab\},a\boxplus b=\{1,ab\}$ and hence
$$1\boxplus(1\boxplus a)=1\boxplus\{1,a,b,ab\}=\{1,a,b,ab\}$$
Meanwhile, $a\boxplus (1\boxplus1)=a\boxplus\{a,b\}=\{1,ab\}$, which is impossible. So there are no valid hyperfields in this case.\\
\textbf{Case 9} If $1\boxplus1=\{1,a\}$, then $1\boxplus a=\{1,a\}$ and
$$1\boxplus ab=\{1,a,b,ab,0\},1\boxplus b=\{1,a,b,ab\}$$
Therefore one can get the following hyperfield.
$$\begin{tabular}{|c|c|c|c|c|c|}
 \hline
  $\boxplus$ & $0$ & $1$ & $a$ & $b$ & $ab$\\
  \hline
  $0$ & $\{0\}$ & $\{1\}$ & $\{a\}$ & $\{b\}$ & $\{ab\}$\\
  \hline
  $1$ & $\{1\}$ & $\{1,a\}$ & $\{1,a\}$ & $\{1,a,b,ab\}$ & $\{1,a,b,ab,0\}$\\
  \hline
  $a$ & $\{a\}$  & $\{1,a\}$ & $\{1,a\}$ & $\{1,a,b,ab,0\}$ & $\{1,a,b,ab\}$\\
  \hline
  $b$ & $\{b\}$ & $\{1,a,b,ab\}$ & $\{1,a,b,ab,0\}$ & $\{b,ab\}$ & $\{b,ab\}$\\
  \hline
  $ab$ & $\{ab\}$ & $\{1,a,b,ab,0\}$ & $\{1,a,b,ab\}$ & $\{b,ab\}$ & $\{b,ab\}$\\
  \hline
\end{tabular}$$
This hyperfield is an extension of $\mathbb{S}^{\Uparrow5}$, denoted by $\mathbb{S}^{\Uparrow5}_\rightarrow$.\\
\textbf{Case 10} If $1\boxplus1=\{1,b\}$, then $1\boxplus b=\{1,b\}$ and
$$1\boxplus ab=\{1,a,b,ab,0\},1\boxplus a=\{1,a,b,ab\}$$
Therefore, one can get the following hyperfield.
$$\begin{tabular}{|c|c|c|c|c|c|}
 \hline
  $\boxplus$ & $0$ & $1$ & $a$ & $b$ & $ab$\\
  \hline
  $0$ & $\{0\}$ & $\{1\}$ & $\{a\}$ & $\{b\}$ & $\{ab\}$\\
  \hline
  $1$ & $\{1\}$ & $\{1,b\}$ & $\{1,a,b,ab\}$ & $\{1,b\}$ & $\{1,a,b,ab,0\}$\\
  \hline
  $a$ & $\{a\}$  & $\{1,a,b,ab\}$ & $\{a,ab\}$ & $\{1,a,b,ab,0\}$ & $\{a,ab\}$\\
  \hline
  $b$ & $\{b\}$ & $\{1,b\}$ & $\{1,a,b,ab,0\}$ & $\{1,b\}$ & $\{1,a,b,ab\}$\\
  \hline
  $ab$ & $\{ab\}$ & $\{1,a,b,ab,0\}$ & $\{a,ab\}$ & $\{1,a,b,ab\}$ & $\{a,ab\}$\\
  \hline
\end{tabular}$$
This hyperfield is isomorphic to $\mathbb{S}^{\Uparrow5}_\rightarrow$.\\
\textbf{Case 11} If $1\boxplus1=\{1,a,b\}$, then $\{1,ab\}\subseteq1\boxplus ab$.\\
\textbf{Case 11.1} When $1\boxplus ab=\{1,ab,0\}$, one can obtain that $1\boxplus ab=\{1,a,ab\}$. In this case, we have
$$(1\boxplus ab)\boxplus ab=\{1,ab,0\}\boxplus ab=\{1,a,ab,0\}$$
while
$$(ab\boxplus ab)\boxplus 1=\{1,a,ab\}\boxplus 1=\{1,a,b,ab,0\}$$
It is possible because we need the hyperaddition $\boxplus$ to be associative.\\
\textbf{Case 11.2} When $1\boxplus ab=\{1,a,b,ab,0\}$, we have the following hyperfield.
$$\begin{tabular}{|c|c|c|c|c|c|}
 \hline
  $\boxplus$ & $0$ & $1$ & $a$ & $b$ & $ab$\\
  \hline
  $0$ & $\{0\}$ & $\{1\}$ & $\{a\}$ & $\{b\}$ & $\{ab\}$\\
  \hline
  $1$ & $\{1\}$ & $\{1,a,b\}$ & $\{1,a,b,ab\}$ & $\{1,a,b,ab\}$ & $\{1,a,b,ab,0\}$\\
  \hline
  $a$ & $\{a\}$  & $\{1,a,b,ab\}$ & $\{1,a,ab\}$ & $\{1,a,b,ab,0\}$ & $\{1,a,b,ab\}$\\
  \hline
  $b$ & $\{b\}$ & $\{1,a,b,ab\}$ & $\{1,a,b,ab,0\}$ & $\{1,b,ab\}$ & $\{1,a,b,ab\}$\\
  \hline
  $ab$ & $\{ab\}$ & $\{1,a,b,ab,0\}$ & $\{1,a,b,ab\}$ & $\{1,a,b,ab\}$ & $\{1,a,ab\}$\\
  \hline
\end{tabular}$$
It is an extension of $\mathbb{S}^{\Uparrow5}_\rightarrow$ and we denote it by $\mathbb{S}^{\Uparrow5}_{\rightarrow^2}$.\\
\textbf{Case 12} If $1\boxplus1=\{1,a,ab\}$, then $1\boxplus a=\{1,a\}$ and
$$1\boxplus ab=\{1,a,b,ab,0\},\quad 1\boxplus b=\{1,a,b,ab\}$$
Therefore, one can see that
$$(1\boxplus a)\boxplus 1=\{1,a\}\boxplus 1=\{1,a,ab,0\}$$
while
$$(1\boxplus 1)\boxplus a=\{1,a,ab\}\boxplus 1=\{1,a,b,ab,0\}$$
It is possible because we need the hyperaddition $\boxplus$ to be associative.\\
\textbf{Case 13} If $1\boxplus1=\{1,b,ab\}$, then $1\boxplus b=\{1,b\}$ and
$$1\boxplus ab=\{1,a,b,ab,0\},\quad 1\boxplus a=\{1,a,b,ab\}$$
Therefore, one can see that
$$(1\boxplus b)\boxplus 1=\{1,b\}\boxplus 1=\{1,b,ab,0\}$$
while
$$(1\boxplus 1)\boxplus b=\{1,b,ab\}\boxplus 1=\{1,a,b,ab,0\}$$
It is possible because we need the hyperaddition $\boxplus$ to be associative.\\
\textbf{Case 14} If $1\boxplus1=\{a,b,ab\}$, then $1\boxplus ab=\{a,b,0\}$. We consider set $1\boxplus a$.\\
\textbf{Case 14.1} When $1\boxplus a=\{b,ab\}$ and $1\boxplus b=\{a,ab\}$, one can see
$$ab\in a\boxplus(1\boxplus1)=a\boxplus\{a,b,ab\}=a\boxplus a\cup a\boxplus\{b,ab\}=\{1,b,ab\}\cup a\boxplus\{b,ab\}$$
but $1\boxplus(1\boxplus a)=1\boxplus\{b,ab\}=\{1,a,b,0\}$, a contradiction.\\
\textbf{Case 14.2} When $1\boxplus a=1\boxplus b=\{1,a,b,ab\}$, one can see
$$\begin{tabular}{|c|c|c|c|c|c|}
 \hline
  $\boxplus$ & $0$ & $1$ & $a$ & $b$ & $ab$\\
  \hline
  $0$ & $\{0\}$ & $\{1\}$ & $\{a\}$ & $\{b\}$ & $\{ab\}$\\
  \hline
  $1$ & $\{1\}$ & $\{a,b,ab\}$ & $\{1,a,b,ab\}$ & $\{1,a,b,ab\}$ & $\{a,b,0\}$\\
  \hline
  $a$ & $\{a\}$  & $\{1,a,b,ab\}$ & $\{1,b,ab\}$ & $\{1,ab,0\}$ & $\{1,a,b,ab\}$\\
  \hline
  $b$ & $\{b\}$ & $\{1,a,b,ab\}$ & $\{1,ab,0\}$ & $\{a,b,ab\}$ & $\{1,a,b,ab\}$\\
  \hline
  $ab$ & $\{ab\}$ & $\{a,b,0\}$ & $\{1,a,b,ab\}$ & $\{1,a,b,ab\}$ & $\{1,b,ab\}$\\
  \hline
\end{tabular}$$
From my perspective, this hyperfield is a special one since it is not an extension of any hyperfield we mentioned before. Here I denote it by $\mathbb{X}$.\\
\textbf{Case 15} If $1\boxplus1=\{1,a,b,ab\}$, then $1\boxplus ab=\{1,a,b,ab,0\}$, we consider set $1\boxplus a$.\\
\textbf{Case 15.1} When $1\boxplus a=\{b,ab\}$ and $1\boxplus b=\{a,ab\}$, one can see
$$ab\in a\boxplus(1\boxplus1)=a\boxplus\{a,b,ab\}=a\boxplus a\cup a\boxplus\{b,ab\}=\{1,b,ab\}\cup a\boxplus\{b,ab\}$$
but $1\boxplus(1\boxplus a)=1\boxplus\{b,ab\}=\{1,a,b,0\}$, a contradiction.

\noindent\textbf{Case 15.2} When $1\boxplus a=1\boxplus b=\{1,a,b,ab\}$, we have hyperfield
$$\begin{tabular}{|c|c|c|c|c|c|}
 \hline
  $\boxplus$ & $0$ & $1$ & $a$ & $b$ & $ab$\\
  \hline
  $0$ & $\{0\}$ & $\{1\}$ & $\{a\}$ & $\{b\}$ & $\{ab\}$\\
  \hline
  $1$ & $\{1\}$ & $\{1,a,b,ab\}$ & $\{1,a,b,ab\}$ & $\{1,a,b,ab\}$ & $\{1,a,b,ab,0\}$\\
  \hline
  $a$ & $\{a\}$  & $\{1,a,b,ab\}$ & $\{1,a,b,ab\}$ & $\{1,a,b,ab,0\}$ & $\{1,a,b,ab\}$\\
  \hline
  $b$ & $\{b\}$ & $\{1,a,b,ab\}$ & $\{1,a,b,ab,0\}$ & $\{1,a,b,ab\}$ & $\{1,a,b,ab\}$\\
  \hline
  $ab$ & $\{ab\}$ & $\{1,a,b,ab,0\}$ & $\{1,a,b,ab\}$ & $\{1,a,b,ab\}$ & $\{1,a,b,ab\}$\\
  \hline
\end{tabular}$$
We denote this hyperfield by $\mathbb{W}^{\Uparrow5}_{\rightarrow}$.
\end{proof}
\end{ppst}
\begin{rem}
One can check the following commutative diagram.
$$
\xymatrix{
\mathbb{S}\ar[rr]\ar[d] &&\mathbb{S}^{\Uparrow5}\ar[rr]\ar[d] &&\mathbb{S}^{\Uparrow5}_{\rightarrow}\ar[rr] &&\mathbb{S}^{\Uparrow5}_{\rightarrow^2} \ar[dll] \\
\mathbb{W}\ar[rr] && \mathbb{W}^{\Uparrow5}\ar[rr] && \mathbb{W}^{\Uparrow5}_{\rightarrow} && \mathbb{X}\ar[ll]
}
$$
In addition, by the two isomorphisms
$$f:(1,a,b,c,0)\mapsto(1,b,a,c,0),\quad g:(1,a,b,c)\mapsto(1,c,b,a)$$
where $c=ab$, one can see that each hyperfield with $a=-1$ or $b=-1$ of underlying monoid $C_{2,2}\cup\{0\}$ is isomorphic to one of hyperfields in Proposition \ref{F5,3}.
\end{rem}

\begin{ppst}\label{F5,4}
There are five different hyperfields with $1=-1$ of underlying monoid $C_{2,2}\cup\{0\}$.

\begin{proof}
Here we divide cases by set $1\boxplus1$. Notice that in the monoid $C_{2,2}\cup\{0\}$, one can permutate $a,b,ab$ without changing the structure of $C_{2,2}\cup\{0\}$. So it is sufficient for us to only consider the following cases.\\
\textbf{Case 1} If $1\boxplus1=\{0\}$, then $\{1,a\}\cap1\boxplus a=\varnothing$ and hence $1\boxplus a=\{b,ab\}$. In this case, we have $1\boxplus b=\{a,ab\}$ and then
$$b\in1\boxplus(1\boxplus a)=1\boxplus\{b,ab\}$$
In the meanwhile, $a\boxplus(1\boxplus1)=a\boxplus 0=\{a\}$, a contradiction.\\
\textbf{Case 2} If $1\boxplus1=\{1,0\}$, then $\{1,a\}\cap1\boxplus a=\varnothing$ and hence $1\boxplus a=\{b,ab\}$. In this case, we have the following hyperfield.
$$\begin{tabular}{|c|c|c|c|c|c|}
 \hline
  $\boxplus$ & $0$ & $1$ & $a$ & $b$ & $ab$\\
  \hline
  $0$ & $\{0\}$ & $\{1\}$ & $\{a\}$ & $\{b\}$ & $\{ab\}$\\
  \hline
  $1$ & $\{1\}$ & $\{1,0\}$ & $\{b,ab\}$ & $\{a,ab\}$ & $\{a,b\}$\\
  \hline
  $a$ & $\{a\}$  & $\{b,ab\}$ & $\{a,0\}$ & $\{1,ab\}$ & $\{1,b\}$\\
  \hline
  $b$ & $\{b\}$ & $\{a,ab\}$ & $\{1,ab\}$ & $\{b,0\}$ & $\{1,a\}$\\
  \hline
  $ab$ & $\{ab\}$ & $\{a,b\}$ & $\{1,b\}$ & $\{1,a\}$ & $\{ab,0\}$\\
  \hline
\end{tabular}$$
It is an extension of $\mathbb{K}$, denoted $\mathbb{K}^{\Uparrow5}$.\\
\textbf{Case 3} If $1\boxplus1=\{a,0\}$, then $\{1,b\}\cap1\boxplus b=\varnothing$ and hence $1\boxplus b=\{a,ab\}$. In this case, one can see $1\boxplus a=\{1,a,b,ab\}$. Therefore
$$1\boxplus(1\boxplus a)=1\boxplus\{1,a,b,ab\}=\{1,a,b,ab,0\}$$
while $a\boxplus(1\boxplus1)=a\boxplus\{a,0\}=\{1,a,0\}$, a contradiction.\\
\textbf{Case 4} If $1\boxplus1=\{1,a,0\}$, then $\{1,b\}\cap1\boxplus b=\varnothing$ and hence $1\boxplus b=\{a,ab\}$. In this case, one can see $1\boxplus a=\{1,a,b,ab\}$ since
$$\{a,ab\}\subseteq1\boxplus b\Leftrightarrow \{b,ab\}\subseteq1\boxplus a\Leftrightarrow \{a,b\}\subseteq1\boxplus ab$$
Therefore, we have $ab\boxplus (1\boxplus b)=ab\boxplus \{a,ab\}=\{1,b,ab,0\}$ but
$$1\boxplus (b\boxplus ab)=1\boxplus\{1,a,b,ab\}=\{0,1,a,b,ab\}$$
a contradiction.\\
\textbf{Case 5} If $1\boxplus1=\{a,b,0\}$, then $\{1,ab\}\cap 1\boxplus ab=\varnothing$ and thus $1\boxplus ab=\{1,ab\}$. In this case, we have $1\boxplus a=\{1,a,b,ab\}$ and thus $(1\boxplus a)\boxplus1=\{1,a,b,ab,0\}$. However, one can see $(1\boxplus 1)\boxplus a=\{a,b,0\}\boxplus a=\{1,ab,0\}\cup\{1,ab\}=\{1,a,ab,0\}$, a contradiction.\\
\textbf{Case 6} If $1\boxplus1=\{1,a,b,0\}$, then $\{1,ab\}\cap 1\boxplus ab=\varnothing$ and thus $1\boxplus ab=\{1,ab\}$. In this case, we have $1\boxplus a=1\boxplus b=\{1,a,b,ab\}$ and then can obtain
$$(1\boxplus ab)\boxplus b\neq(1\boxplus b)\boxplus ab$$
So it does not exist a hyperfield in this case.\\
\textbf{Case 7} If $1\boxplus1=\{a,b,ab,0\}$, then $\{1,a\}\subseteq1\boxplus a$.\\
\textbf{Case 7.1} When $1\boxplus a=\{1,a\}$, one can see $1\boxplus b=\{1,b\}$ and $1\boxplus ab=\{1,ab\}$. In this case, one can get the following hyperfield.
$$\begin{tabular}{|c|c|c|c|c|c|}
 \hline
  $\boxplus$ & $0$ & $1$ & $a$ & $b$ & $ab$\\
  \hline
  $0$ & $\{0\}$ & $\{1\}$ & $\{a\}$ & $\{b\}$ & $\{ab\}$\\
  \hline
  $1$ & $\{1\}$ & $\{a,b,ab,0\}$ & $\{1,a\}$ & $\{1,b\}$ & $\{1,ab\}$\\
  \hline
  $a$ & $\{a\}$  & $\{1,a\}$ & $\{1,b,ab,0\}$ & $\{1,ab\}$ & $\{a,ab\}$\\
  \hline
  $b$ & $\{b\}$ & $\{1,b\}$ & $\{1,ab\}$ & $\{1,a,ab,0\}$ & $\{b,ab\}$\\
  \hline
  $ab$ & $\{ab\}$ & $\{1,ab\}$ & $\{a,ab\}$ & $\{b,ab\}$ & $\{1,a,b,0\}$\\
  \hline
\end{tabular}$$
This hyperfield is an extension of $(\mathbb{F}_2)^{\uparrow3}$, denoted by $(\mathbb{F}_2)^{\Uparrow5}$. As one can see, it is also the same sort of hyperfields as $\mathbb{M}$ described in \cite{6}. In that paper, Massouros proves that it is not isomorphic to any quotient of fields.

\noindent\textbf{Case 7.2} When $1\boxplus a=\{1,a,b,ab\}$, one can see $1\boxplus b=1\boxplus ab=\{1,a,b,ab\}$. In this case, one can get the following hyperfield.
$$\begin{tabular}{|c|c|c|c|c|c|}
 \hline
  $\boxplus$ & $0$ & $1$ & $a$ & $b$ & $ab$\\
  \hline
  $0$ & $\{0\}$ & $\{1\}$ & $\{a\}$ & $\{b\}$ & $\{ab\}$\\
  \hline
  $1$ & $\{1\}$ & $\{a,b,ab,0\}$ & $\{1,a,b,ab\}$ & $\{1,a,b,ab\}$ & $\{1,a,b,ab\}$\\
  \hline
  $a$ & $\{a\}$  & $\{1,a,b,ab\}$ & $\{1,b,ab,0\}$ & $\{1,a,b,ab\}$ & $\{1,a,b,ab\}$\\
  \hline
  $b$ & $\{b\}$ & $\{1,a,b,ab\}$ & $\{1,a,b,ab\}$ & $\{1,a,ab,0\}$ & $\{1,a,b,ab\}$\\
  \hline
  $ab$ & $\{ab\}$ & $\{1,a,b,ab\}$ & $\{1,a,b,ab\}$ & $\{1,a,b,ab\}$ & $\{1,a,b,0\}$\\
  \hline
\end{tabular}$$
This is an extension of $(\mathbb{F}_2)^{\Uparrow5}$, denoted by $(\mathbb{F}_2)^{\Uparrow5}_\rightarrow$.\\
\textbf{Case 8} If $1\boxplus1=\{1,a,b,ab,0\}$, then $\{1,a\}\subseteq1\boxplus a$.\\
\textbf{Case 8.1} When $1\boxplus a=\{1,a\}$, one can see $1\boxplus b=\{1,b\}$ and $1\boxplus ab=\{1,ab\}$. In this case, one can get the following hyperfield.
$$\begin{tabular}{|c|c|c|c|c|c|}
 \hline
  $\boxplus$ & $0$ & $1$ & $a$ & $b$ & $ab$\\
  \hline
  $0$ & $\{0\}$ & $\{1\}$ & $\{a\}$ & $\{b\}$ & $\{ab\}$\\
  \hline
  $1$ & $\{1\}$ & $\{1,a,b,ab,0\}$ & $\{1,a\}$ & $\{1,b\}$ & $\{1,ab\}$\\
  \hline
  $a$ & $\{a\}$  & $\{1,a\}$ & $\{1,a,b,ab,0\}$ & $\{1,ab\}$ & $\{a,ab\}$\\
  \hline
  $b$ & $\{b\}$ & $\{1,b\}$ & $\{1,ab\}$ & $\{1,a,b,ab,0\}$ & $\{b,ab\}$\\
  \hline
  $ab$ & $\{ab\}$ & $\{1,ab\}$ & $\{a,ab\}$ & $\{b,ab\}$ & $\{1,a,b,ab,0\}$\\
  \hline
\end{tabular}$$
This is an extension of $\mathbb{K}^{\uparrow3}$, denoted by $(\mathbb{K})^{\Rsh5}$.\\
\textbf{Case 8.2} When $1\boxplus a=\{1,a,b,ab\}$, one can see $1\boxplus b=1\boxplus ab=\{1,a,b,ab\}$. In this case, one can get the following hyperfield.
$$\begin{tabular}{|c|c|c|c|c|c|}
 \hline
  $\boxplus$ & $0$ & $1$ & $a$ & $b$ & $ab$\\
  \hline
  $0$ & $\{0\}$ & $\{1\}$ & $\{a\}$ & $\{b\}$ & $\{ab\}$\\
  \hline
  $1$ & $\{1\}$ & $\{1,a,b,ab,0\}$ & $\{1,a,b,ab\}$ & $\{1,a,b,ab\}$ & $\{1,a,b,ab\}$\\
  \hline
  $a$ & $\{a\}$  & $\{1,a,b,ab\}$ & $\{1,a,b,ab,0\}$ & $\{1,a,b,ab\}$ & $\{1,a,b,ab\}$\\
  \hline
  $b$ & $\{b\}$ & $\{1,a,b,ab\}$ & $\{1,a,b,ab\}$ & $\{1,a,b,ab,0\}$ & $\{1,a,b,ab\}$\\
  \hline
  $ab$ & $\{ab\}$ & $\{1,a,b,ab\}$ & $\{1,a,b,ab\}$ & $\{1,a,b,ab\}$ & $\{1,a,b,ab,0\}$\\
  \hline
\end{tabular}$$
This is an extension of $\mathbb{K}^{\Uparrow5}$, denoted by $\mathbb{K}^{\Uparrow5}$.
\end{proof}
\end{ppst}
\begin{rem}
One can check the following commutative diagram.
$$
\xymatrix{
\mathbb{F}_2\ar[r]\ar[d] &(\mathbb{F}_2)^{\uparrow3}\ar[rrr]\ar[dd]&  &&(\mathbb{F}_2)^{\Uparrow5}\ar[r]&(\mathbb{F}_2)^{\Uparrow5}_{\rightarrow} \ar[d] \\
\mathbb{K}\ar[rr]\ar[dr] && \mathbb{K}^{\Uparrow5}\ar[rrr] &&& \mathbb{K}^{\Uparrow5}_\rightarrow \\
&\mathbb{K}^{\uparrow3}\ar[rr]&&(\mathbb{K})^{\Rsh5}\ar[urr]
}
$$
\end{rem}
\begin{thm}
There are, up to hyperfield isomorphisms, exactly 27 hyperfields of order 5.

\begin{proof}
According to Proposition \ref{F5,1}, Proposition \ref{F5,2}, Proposition \ref{F5,3} and Proposition \ref{F5,4}, we have $9+7+6+5=27$ different hyperfields of order 5.
\end{proof}
\end{thm}
\begin{rem}
For more information and a deeper point of view, readers can refer to Mansour Eyvazi's work \cite{2}. He gives an algorithm to compute finite hyperfields in that work.
\end{rem}


\end{document}